\documentclass[11pt]{article}

\usepackage{algorithm,algpseudocode}
\usepackage[utf8]{inputenc}
\usepackage{fullpage}
\usepackage{amsmath,amsthm,amssymb}
\usepackage{mathabx}
\usepackage{bm}
\usepackage{graphicx}
\usepackage{float}
\usepackage[margin=1in]{geometry}
\usepackage{tikz}
\usepackage{hyperref}
\usepackage{mathtools}
\usepackage{bbm}
\usepackage{booktabs}
\usepackage{soul}
\usepackage{authblk}

\newcommand{\mr}{\mathrm}
\newcommand{\mc}{\mathcal}
\newcommand{\mb}{\mathbf}

\newcommand{\pp}{\partial}

\theoremstyle{remark}

\title{Element learning: a systematic approach of accelerating finite element-type methods via machine learning, with applications to radiative transfer}
\author[1]{Shukai Du\thanks{Email: sdu49@wisc.edu}}
\author[1,2]{Samuel N. Stechmann}
\affil[1]{Department of Mathematics, University of Wisconsin–Madison}
\affil[2]{Department of Atmospheric and Oceanic Sciences, University of Wisconsin–Madison}

\date{Revised May 31, 2024}

\begin{document}

\maketitle

\begin{abstract}
In this paper, we propose a systematic approach for accelerating finite element-type methods by machine learning for the numerical solution of partial differential equations (PDEs).  
The main idea is to use a neural network to learn the solution map of the PDEs and to do so in an element-wise fashion.
This map takes input of the element geometry and the PDE's parameters on that element, and gives output of two operators: (1) the {\sl in2out} operator for inter-element communication, and (2) the {\sl in2sol}  operator (Green's function) for element-wise  solution recovery.
A significant advantage of this approach is that, once trained, this network can be used for the numerical solution of the PDE for any domain geometry and any parameter distribution without retraining. 
Also, the training is significantly simpler since it is done on the element level instead on the entire domain.
We call this approach element learning. This method is closely related to hybridizable discontinuous Galerkin (HDG) methods in the sense that the local solvers of HDG are replaced by machine learning approaches. 
%As a result, all HDG methods, along with hybridized mixed and hybridized continuous Galerkin methods can be potentially accelerated by the proposed element learning approach.

Numerical tests are presented for an example PDE, the radiative transfer or radiation transport equation, in a variety of scenarios with idealized or realistic cloud fields, with smooth or sharp gradient in the cloud boundary transition.
%including the radiative distribution calculation for an idealized cloud with smooth or sharp cloud boundary transition, and a more realistic cloud field from the I3RC test cases. 
Under a fixed accuracy level of $10^{-3}$ in the relative $L^2$ error, and polynomial degree $p=6$ in each element, we observe an approximately 5 to 10 times speed-up by element learning compared to a classical finite element-type method.
%in a relatively high-$p$ regime of spectral element.
%Remaining considerations include exploring the greater speed-up that can be expected if even larger $p$ values are used.

\vspace{0.5cm}\noindent
\end{abstract}

\section{Introduction}
\label{sec:intro}

\subsection{Background and motivation}
In the past decade, (artificial) neural networks and machine learning tools have surfaced as game-changing technologies across numerous fields, resolving an array of challenging problems. Examples include image recognition \cite{lecun2015deep,krizhevsky2017imagenet}, playing the game go \cite{silver2016mastering}, protein folding \cite{jumper2021highly}, and large language models such as GPT3 \cite{brown2020language}. 

Given these impressive results, it is reasonable to envision the potential of neural networks (NNs)
for the numerical solution of partial differential equations (PDEs) or other scientific computing problems. There has been significant work in this direction \cite{brunton2016discovering,pathak2018model,HeLiXuZh:2020,fan2019multiscale,vlachas2020backpropagation,qi2020using,wu2020enforcing,kochkov2021machine,khoo2021solving,chen2021bamcafe,chen2021reduced,chen2021solving,gottwald2021combining,hu2021structure,TaRa:2021,zepeda2021deep,ding2022coupling,dehoop2022cost,engquist2022generalized,kaneda2022deep,schneider2022ensemble,yu2022data,darcy2023one,gregory2023deep,horenko2023cheap,linot2023stabilized,park2023near,batlle2024kernel,bruna2024neural}. As  examples, we mention two major groups  -- (1) neural networks as function approximators, (2) neural networks as operator approximators.

The first type of these methods uses neural networks to approximate the solution of PDEs. Examples include Physics-Informed Neural Networks (PINNs) \cite{raissi2019physics,mishra2021physics,lu2022solving,cuomo2022scientific} and Deep Ritz methods \cite{yu2018deep}. These methods have demonstrated promising results, especially in the realm of high-dimensional problems or inverse problems, since they seem to bypass the notorious `curse of dimensionality', potentially attributable to the neural network's efficient approximation capabilities in high-dimensional spaces \cite{han2017deep,han2018solving,MoYaDu:2021}.

While the potential advantages are promising, it remains uncertain if the aforementioned methods have a real advantage over traditional algorithms, such as finite element-type methods, when addressing many classic PDEs. A crucial factor contributing to this uncertainty stems from the complex optimization problems and error landscapes—often non-convex—introduced by neural networks due to their hidden layers and non-linear activation functions. 
%This complexity becomes particularly notable when considering that many of these PDEs are linear in nature, including Maxwell's equations, radiative transfer, linear elasticity, and acoustic waves. Contrarily, neural networks are inherently non-linear representations. This disparity can complicate the optimization process, resulting in slow training or inadequate error reduction throughout the training phase,  due to the intricate error landscape of the optimization function.

The second type of these methods employs neural networks to approximate operators, with methods like Neural Operator \cite{kovachki2023neural,molinaro2023neural} and DeepOnet \cite{lu2021learning} serving as leading examples. A primary advantage of these methods hinges on their speed -- once trained, solving a problem is simply a forward propagation of the network. 
For instance, in \cite{pathak2022fourcastnet}, it was reported that the neural operator based method can be {orders of magnitude} faster than classical finite volume/difference based methods in weather simulation.
{This suggests great potential of using machine learning to accelerate computation.}

However, this type of approach is essentially data-driven and its performance can depend on the training samples and how well it generalizes. Consequently, the reliability of the result may falter during extreme events, as indicated in \cite{pathak2022fourcastnet}.
In addition to this, {\color{black} 
it is common for these approaches to be constrained} by the geometry of the domain and boundary conditions. 
For instance, the Fourier Neural Operator method that relies on Fourier expansion for solution representation is specifically tailored for rectangular/cubic domains with periodic boundary conditions \cite{kovachki2023neural}. 
{\color{black}However, there are recent attempts to extend operator learning to more complex geometries \cite{serrano2024operator,li2024geometry}.

}

%Despite the PDEs might be linear, the map from the parameter space to the solution space is usually non-linear. Based on this observation

On the other hand, for many years, traditional methods like finite element-type methods  have played a critical role in the advancement of various scientific and engineering disciplines. These methodologies, honed and well-understood over more than five decades, have given us a plethora of reliable and robust techniques successfully employed across an array of problems. Examples include, but are definitely not limited to, conforming and non-conforming finite element \cite{BrSc:2008,CrRa:1973}, mixed finite element \cite{RaTh:1977,Ne:1980}, discontinuous Galerkin \cite{CoSh:2001,ArBrCo:2001}, along with effective techniques such as slope limiter \cite{Le:1992,CoSh:2001}, multi-scale finite elements \cite{efendiev2009multiscale}, finite element exterior calculus \cite{ArFaWi:2006,GiHoZh:2017}, $hp$-adaptivity \cite{BaSu:1994,HoSu:2005} to tailor these methods to different application scenarios \cite{Mo:2003,duben2012discontinuous,heath2012discontinuous,cotter2014finite,du2023fast,du2023inverse,MaKeMo:2016,TaTrIs:1997,GaRj:2018}.

Despite the undeniable utility of these classical approaches, they also come with constraints. A primary limitation is on their speed, especially when facing high-dimensional problems, such as those found in radiative transfer, or scenarios when fast forward solvers are necessary, such as in inverse problems. Furthermore, they may lack efficiency in dealing with data or geometric structures characterized by multi-scale variations.
Considering these factors, it is desirable to devise ways of reducing the computational cost of these classical approaches.

%\st{In conclusion, many current approaches employ neural networks either as an approximator for the PDEs' solution or as a data-driven method for approximating solution operators, without incorporating much existing techniques from traditional numerical methods.}

In this paper, we propose a novel approach (which we call \textit{element learning}), {aimed at} accelerating finite element-type methods via machine learning. {This effort seeks to leverage the extensive knowledge gained from decades of development in classical finite element methods}, with the objective of using the modern computational power of machine learning to enhance the computational speed.

%\st{Let us clarify our targets. We are not trying to tackle extremely high-dimensional problems (e.g., with $10^3$ dimensions or more), as neural networks appear to have a significant advantage over finite element-type methods in these scenarios. Instead, our emphasis is on moderately-high dimensional problems such as radiative transfer (with 5 dimensions, 3 spatial and 2 angular), or on inverse problems that require fast forward solvers. In these contexts, the efficiency of traditional solvers might not fulfill application requirements, whereas neural networks may falter in delivering satisfactory outcomes, which is often due to intricate domain geometries, multiscale variations in coefficients or input data, or the current lack of thorough theoretical understanding of neural networks, rendering them not satisfactorily reliable.}

\subsection{Main idea}

To describe the main ideas of the methods here, we will utilize the contexts of both a general class of partial differential equations (PDEs) and also the specific example of the radiative transfer or radiation transport equation.
To proceed with the discussion, we note that a wide variety of PDEs, including radiative transfer, can be written in the following abstract formulation:
\begin{align}\label{eq:PDE}
\mathcal{L}(\sigma)[u] = f \quad \text{in} \quad \Omega.
\end{align}
In this representation, $\mathcal{L}(\sigma)$ represents a combination of differential and integral operators, which also depend on the coefficients $\sigma$. We are in pursuit of a solution denoted by $u$, while $f$ is given as the data, comprising both the external force term and the prescribed boundary conditions.
For example, in the context of atmospheric radiation, $\sigma$ represents the optical properties of the atmosphere (extinction and scattering coefficients), $u$ denotes the radiant intensity, and $f$ symbolizes the inflow and black body radiation.

%{As an example problem for this first formulation and application of element learning}, our research shall centers on radiative transfer equations. 
The specific PDE example here is chosen to be the radiative transfer equation for a number of reasons.
%This equation sees a diverse utilization across various 
Radiative transfer is important in many
fields, including medical imaging \cite{ArSc:2009,Re:2010}, neutron transport \cite{reed1973triangular,lesaint1974finite}, and climate and weather prediction \cite{hogan2017radiation,hogan2018flexible}. Furthermore, the equation offers intriguing mathematical characteristics as it exhibits both hyperbolic properties in weak scattering scenarios and elliptic properties in instances of strong scattering. Finally, the equation is computational challenging because of its high or moderately-high dimensionality so a faster solver is still highly demanded.

% The  solution to the equation \eqref{eq:PDE}  may encounter two primary impediments: (1) the emergence of multi-scale characteristics within the coefficient $\sigma$ or the data $f$, and (2) potential geometric intricacies intrinsic to the domain $\Omega$.

% Finite element-type methods confront these challenges using an old but effective strategy — `divide and conquer'. 
% This involves partitioning the domain into a collection of elements, denoted as $K$, thereby reducing the variations originating from $\sigma$ or $f$ within each element. This partitioning also greatly reduces geometric complexity, with each element $K$ potentially represented as a simple triangle (in $\mathbb R^2$) or a tetrahedron (in $\mathbb R^3$). Consequently, the solution $u$ and the data $\sigma$ and $f$ can be better approximated by polynomials, and solving the PDEs within each element becomes a substantially simpler task than addressing the PDEs across the entire domain.

For finite element-type methods, the domain $\Omega$ is partitioned into a collection of simple geometric objects (e.g. triangles, tetrahedron, rectangles, cubes), and a key aspect for the success of these methods is the effective incorporation of the solution operators on each element to construct the solution for the entire domain. %Essential to this process is the management of communication between elements, which 
The study of this 
has engendered a plethora of diverse finite element methods, including conforming and non-conforming, mixed finite element, and discontinuous Galerkin methods. 
Recently, it was realised that many of these methods can be implemented within the framework of hybridizable discontinuous Galerkin (HDG) methods in a unified way \cite{CoGoLa:2009,DuSa:2019,CoFuSa:2017}.

To elaborate on hybridizable finite element-type methods and their key aspects for element learning, we would need to find the discretization of the following two operators on an element-wise basis:

\begin{itemize}
\item The in2out operator responsible for inter-element communication:
\begin{align}\label{eq:in_out}
    u_{\pp K}^\mr{in}\rightarrow u_{\pp K}^\mr{out}.
\end{align}
For radiative transfer, it is the inflow-radiation-to-outflow-radiation operator.
For elliptic-type equations, this equates to the Dirichlet-to-Neumann operator (the DtN map). 
\item The in2sol operator responsible for recovering the solution (or desired properties of the solution) on each element: 
\begin{align}\label{eq:in_sol}
    u_{\pp K}^\mr{in}\rightarrow u_K.
\end{align}
For radiative transfer, it is the inflow-radiation-to-solution operator.
For elliptic-type equations, this is the Dirichlet-to-solution operator (i.e., the Green's function). Note that we can also replace $u_K$ by its partial statistics if desired. For instance, we can replace the radiative intensity $u_K$ by its mean intensity or radiative heating rate, if desired.
\end{itemize}

\def\inout(#1,#2,#3){
\draw (#1,#2) -- ++(#3,0) -- ++(0,#3) -- ++(-#3,0) -- cycle;
\draw[->] (#1+0.3*#3,#2) -- ++(0,0.2*#3);
\draw[->] (#1+0.6*#3,#2) -- ++(0,0.2*#3);
\draw[->,double,double distance=1.5pt] (#1+0.3*#3,#2) -- ++(0,-0.2*#3);
\draw[->,double,double distance=1.5pt] (#1+0.6*#3,#2) -- ++(0,-0.2*#3);

\draw[->] (#1+0.3*#3,#2+#3) -- ++(0,-0.2*#3);
\draw[->] (#1+0.6*#3,#2+#3) -- ++(0,-0.2*#3);
\draw[->,double,double distance=1.5pt] (#1+0.3*#3,#2+#3) -- ++(0,0.2*#3);
\draw[->,double,double distance=1.5pt] (#1+0.6*#3,#2+#3) -- ++(0,0.2*#3);

\draw[->] (#1+#3,#2+0.3*#3) -- ++(-0.2*#3,0);
\draw[->] (#1+#3,#2+0.6*#3) -- ++(-0.2*#3,0);
\draw[->,double,double distance=1.5pt] (#1+#3,#2+0.3*#3) -- ++(0.2*#3,0);
\draw[->,double,double distance=1.5pt] (#1+#3,#2+0.6*#3) -- ++(0.2*#3,0);

\draw[->] (#1,#2+0.3*#3) -- ++(0.2*#3,0);
\draw[->] (#1,#2+0.6*#3) -- ++(0.2*#3,0);
\draw[->,double,double distance=1.5pt] (#1,#2+0.3*#3) -- ++(-0.2*#3,0);
\draw[->,double,double distance=1.5pt] (#1,#2+0.6*#3) -- ++(-0.2*#3,0);
}

\def\insol(#1,#2,#3){
\draw (#1,#2) -- ++(#3,0) -- ++(0,#3) -- ++(-#3,0) -- cycle;
\draw[->] (#1+0.3*#3,#2) -- ++(0,0.2*#3);
\draw[->] (#1+0.6*#3,#2) -- ++(0,0.2*#3);

\draw[->] (#1+0.3*#3,#2+#3) -- ++(0,-0.2*#3);
\draw[->] (#1+0.6*#3,#2+#3) -- ++(0,-0.2*#3);

\draw[->] (#1+#3,#2+0.3*#3) -- ++(-0.2*#3,0);
\draw[->] (#1+#3,#2+0.6*#3) -- ++(-0.2*#3,0);

\draw[->] (#1,#2+0.3*#3) -- ++(0.2*#3,0);
\draw[->] (#1,#2+0.6*#3) -- ++(0.2*#3,0);

\node at (#1+0.3*#3,#2+0.3*#3) {\scriptsize $\times$};
\node at (#1+0.3*#3,#2+0.6*#3) {\scriptsize $\times$};
\node at (#1+0.6*#3,#2+0.3*#3) {\scriptsize $\times$};
\node at (#1+0.6*#3,#2+0.6*#3) {\scriptsize $\times$};
}

\def\diagexp(#1,#2){
\draw[->] (#1,#2+1.2) -- ++(0.5,0);
\node[align=left,text width=3.5cm] at (#1+2.5,#2+1.2) {inflow};

\draw[->,double,double distance=1.5pt] (#1,#2+0.6) -- ++(0.5,0);
\node[align=left,text width=3.5cm] at (#1+2.5,#2+0.6) {outflow};

\node at (#1+0.25,#2) {$\times$};
\node[align=left,text width=3.5cm] at (#1+2.5,#2) {solution};
}

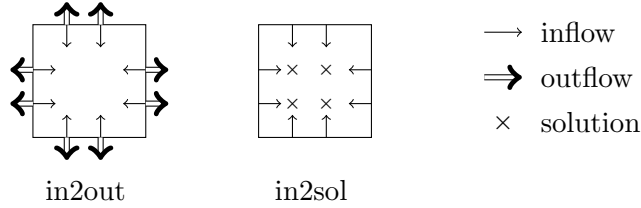
\begin{figure}[H]
\centering
\begin{tikzpicture}
\node at (-4,0) {};

\inout(0,0,1.5)
\node at (0.7,-0.7) {in2out};
\insol(3,0,1.5)
\node at (3.7,-0.7) {in2sol};

\diagexp(6,0.2)
\end{tikzpicture}

\caption{Symbolic representation of the in2out and in2sol operators, from \eqref{eq:in_out} and \eqref{eq:in_sol}, on a square element. When these operators from one element are coupled with the operators from neighboring elements, they can be used to find a global solution (see Figure \ref{fig:hdg_diag}).}
\label{fig:sch_diagr_in2out_in2sol}
\end{figure}

Once these two types of operators or their approximations are obtained, we can recover the solution on the whole domain following the standard HDG routine: (1) use the in2out operators to assemble a global system (involving only the variables on the skeleton of the mesh), (2) solve the global system, (3) use in2sol operator to recover the solution of interest in an element-wise fashion. See the following Figure \ref{fig:hdg_diag} for a visual demonstration.
(Further details will be elaborated in Section \ref{sec:HDG}.)

\def\aarr(#1,#2,#3) {
\draw (#1,#2) circle (2pt);
\draw (#1+0.2*#3,#2) circle (2pt);
\draw (#1+0.5*#3,#2) circle (2pt);
\draw (#1+0.8*#3,#2) circle (2pt);
\draw (#1+#3,#2) circle (2pt);
}
\def\cubea(#1,#2,#3) {
\draw (#1,#2) -- ++(#3,0) -- ++ (0,#3) -- ++(-#3,0) -- cycle;
\aarr(#1,#2,#3)
\aarr(#1,#2+0.2*#3,#3)
\aarr(#1,#2+0.5*#3,#3)
\aarr(#1,#2+0.8*#3,#3)
\aarr(#1,#2+#3,#3)
}
\def\cubeahl(#1,#2,#3) {
\draw (#1,#2) -- ++(#3,0) -- ++ (0,#3) -- ++(-#3,0) -- cycle;
\aarr(#1,#2,#3)
\aarr(#1,#2+#3,#3)
\draw(#1,#2+0.2*#3) circle (2pt);
\draw(#1,#2+0.5*#3) circle (2pt);
\draw(#1,#2+0.8*#3) circle (2pt);
\draw(#1+#3,#2+0.2*#3) circle (2pt);
\draw(#1+#3,#2+0.5*#3) circle (2pt);
\draw(#1+#3,#2+0.8*#3) circle (2pt);
}

\def\aarrcrs(#1,#2,#3) {
\node at (#1,#2) {\small$\times$};
\node at (#1+0.2*#3,#2) {\small$\times$};
\node at (#1+0.5*#3,#2) {\small$\times$};
\node at (#1+0.8*#3,#2) {\small$\times$};
\node at (#1+#3,#2) {\small$\times$};
}
\def\cubeahlcrs(#1,#2,#3) {
\draw (#1,#2) -- ++(#3,0) -- ++ (0,#3) -- ++(-#3,0) -- cycle;
\aarrcrs(#1,#2,#3)
\aarrcrs(#1,#2+#3,#3)
\node at (#1,#2+0.2*#3) {\small $\times$};
\node at (#1,#2+0.5*#3) {\small$\times$};
\node at (#1,#2+0.8*#3) {\small$\times$};
\node at (#1+#3,#2+0.2*#3) {\small$\times$};
\node at (#1+#3,#2+0.5*#3) {\small$\times$};
\node at (#1+#3,#2+0.8*#3) {\small$\times$};
}
\def\cubeacrs(#1,#2,#3) {
\draw (#1,#2) -- ++(#3,0) -- ++ (0,#3) -- ++(-#3,0) -- cycle;
\aarrcrs(#1,#2,#3)
\aarrcrs(#1,#2+0.2*#3,#3)
\aarrcrs(#1,#2+0.5*#3,#3)
\aarrcrs(#1,#2+0.8*#3,#3)
\aarrcrs(#1,#2+#3,#3)
}

\begin{figure}[H]
\centering
\begin{tikzpicture}
%\cubea(0,0,1.5)
%\cubea(0,1.7,1.5)
%\draw[->] (1.7,1.6) -- ++(1.2,0) 
%node [above,midway,sloped] {\scriptsize in2out};
\cubeahl(3.2,1.6,1.5)
\cubeahl(3.2,0.1,1.5)
\node at (3.95,0.85) {\scriptsize in2out};
\node at (3.95,2.35) {\scriptsize in2out};
\draw[->] (5,1.6) -- ++(1.85,0) 
node [above,midway,sloped] {\scriptsize global solve};

\cubeahlcrs(7,0.1,1.5)
\cubeahlcrs(7,1.6,1.5)

\draw[->] (8.7,1.6) -- ++(1.2,0) 
node [above,midway,sloped] {\scriptsize in2sol};

\cubeacrs(10.2,0,1.5)
\cubeacrs(10.2,1.7,1.5)

\end{tikzpicture}
\caption{Solution strategy of both element learning and the standard HDG method. To create a global solver, the local solvers (in2out operators, from Figure \ref{fig:sch_diagr_in2out_in2sol}) from neighboring elements are coupled together. In this schematic diagram, we use `$\circ$' for undetermined DOFs and `$\times$' for determined DOFs (solved solution).}
\label{fig:hdg_diag}
\end{figure}
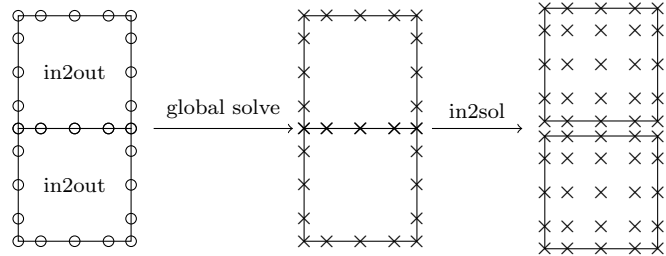

% \begin{figure}[H]
% \centering
% \begin{tikzpicture}
% \cubeahl(0,0,1.5)
% \node at (0.75,0.75) {\scriptsize in2out};
% \cubeahl(0,1.7,1.5)
% \node at (0.75,2.45) {\scriptsize in2out};
% \draw[->] (1.7,1.6) -- ++(1.2,0) 
% node [above,midway,sloped] {\scriptsize coupling};
% \cubeahl(3.2,1.6,1.5)
% \cubeahl(3.2,0.1,1.5)
% \draw[->] (5,1.6) -- ++(1.85,0) 
% node [above,midway,sloped] {\scriptsize global solve};

% \cubeahlcrs(7,0.1,1.5)
% \cubeahlcrs(7,1.6,1.5)

% \draw[->] (8.7,1.6) -- ++(1.2,0) 
% node [above,midway,sloped] {\scriptsize in2sol};

% \cubeacrs(10.2,0,1.5)
% \cubeacrs(10.2,1.7,1.5)

% \end{tikzpicture}
% \caption{A schematic diagram for a standard HDG implementation. Here we use `$\circ$' for undetermined DOFs and `$\times$' for determined DOFs (solved solution).}
% \label{fig:hdg_diag}
% \end{figure}

We now introduce the primary concept of element learning. In this approach, 
%\st{rather than the usually time-consuming approach \cite{woopen2014comparison} of }
instead of directly computing the approximations for the in2out and in2sol operators by inverting a local system on each element, we leverage a neural network to provide the approximations to the in2out and in2sol operators. Note that the neural network will not itself be the approximation of the in2out nor in2sol operator; instead, here, the neural network will provide the in2out and in2sol operators at \textit{outputs}. That is, the neural network approximates the following operators:
\begin{subequations}
\label{eq:nn}
\begin{align}
\label{eq:nn_in_out}
    (\sigma,K) &\xrightarrow{\text{\quad}} \left[u_{\pp K}^\mr{in}\rightarrow u_{\pp K}^\mr{out}\right],\\
\label{eq:nn_in_sol}
    (\sigma,K) &\xrightarrow{\text{\quad}} \left[u_{\pp K}^\mr{in}\rightarrow u_K\right].
\end{align}
\end{subequations}
The above operators take the coefficients $\sigma$ and the element geometry $K$ as inputs. The outputs, in turn, are the in2out and in2sol operators. 
Thus, once this network is trained, it only needs a simple forward propagation of the network to obtain the approximations for the in2out and in2sol operators.
It's important to highlight that although the PDEs might be linear (the in2out and in2sol are linear operators), the operators \eqref{eq:nn} are non-linear. This justifies the use of neural networks with non-linear activation functions and hidden layers.

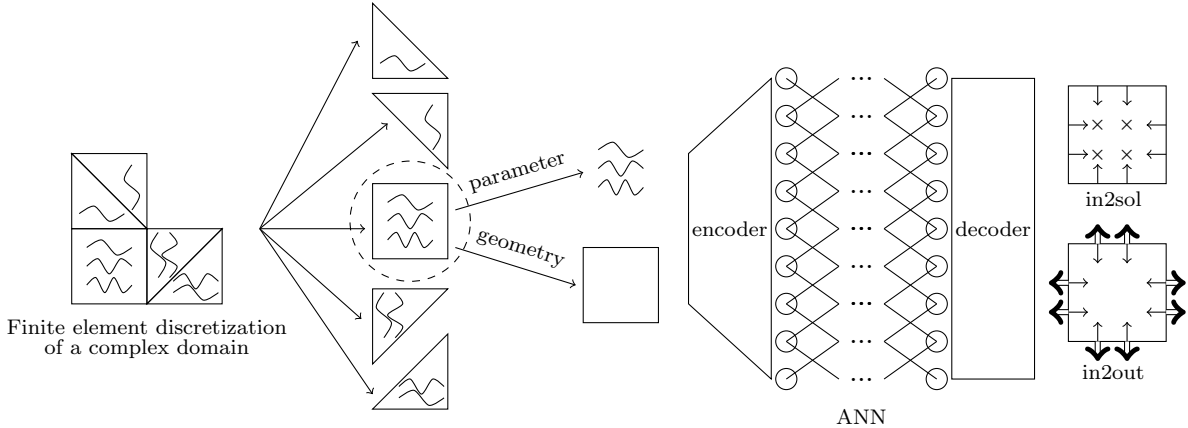
\begin{figure}[H]
\centering
\begin{tikzpicture}
% customized command - curves
\def\crva(#1,#2){
    \draw [black] plot [smooth] coordinates {(#1,#2) (#1+0.2,#2+0.1) (#1+0.4,#2-0.1) (#1+0.6,#2)}
}
\def\crvb(#1,#2){
    \draw [black] plot [smooth] coordinates {(#1,#2) (#1+0.1,#2+0.2) (#1-0.1,#2+0.4) (#1,#2+0.6)}
}
\def\crvbb(#1,#2){
    \draw [black] plot [smooth] coordinates {(#1,#2) (#1+0.1,#2+0.15) (#1-0.1,#2+0.3) (#1+0.1,#2+0.45)
    (#1,#2+0.6)}
}
\def\crvaa(#1,#2){
    \draw [black] plot [smooth] coordinates {(#1,#2) (#1+0.15,#2+0.1) (#1+0.3,#2-0.1) (#1+0.45,#2+0.1)
    (#1+0.6,#2)}
}
\def\crvaaa(#1,#2){
    \draw [black] plot [smooth] coordinates {(#1,#2) (#1+0.12,#2+0.1) (#1+0.24,#2-0.1) (#1+0.36,#2+0.1)
    (#1+0.48,#2-0.1) (#1+0.6,#2)}
}
% customized command - triangle and cubes
\def\triga(#1,#2){\draw  (#1,#2) -- ++(1,0) -- ++(-1,1) -- cycle;
\crva(#1+0.1,#2+0.2);
}
\def\trigb(#1,#2){\draw (#1,#2) -- ++(0,1) -- ++(-1,0) -- cycle;
\crvb(#1-0.2,#2+0.25);
}
\def\trigc(#1,#2){\draw (#1,#2) -- ++(0,1) -- ++(1,0) -- cycle;
\crvbb(#1+0.3,#2+0.35);
\crvb(#1+0.15,#2+0.35);
}
\def\trigd(#1,#2){\draw (#1,#2) -- ++(1,0) -- ++(0,1) -- cycle;
\crvaa(#1+0.35,#2+0.3);
\crva(#1+0.35,#2+0.15);
}
\def\cubea(#1,#2){\draw (#1,#2) -- ++(1,0) -- ++(0,1) -- ++(-1,0) -- cycle;
\crva(#1+0.2,#2+0.75);
\crvaa(#1+0.2,#2+0.5);
\crvaaa(#1+0.2,#2+0.25);
}

% decomposition of the L-shaped domain
\triga(0,0) \trigb(1,0) \cubea(0,-1) \trigc(1,-1) \trigd(1,-1);
\node at (1,-1.3) {\scriptsize Finite element discretization};
\node at (1,-1.6) {\scriptsize of a complex domain};

% decomposition into elements
\draw[->] (2.5,0) -- (3.8,2.5);
\draw[->] (2.5,0) -- (4.2,1.4);
\draw[->] (2.5,0) -- (3.9,0);
\draw[->] (2.5,0) -- (3.85,-1.2);
\draw[->] (2.5,0) -- (4,-2.2);

% the collection of elements
\triga(4,2) \trigb(5,0.8);
\cubea(4,-0.4) 
\trigc(4,-1.8)
\trigd(4,-2.4)

% focus on one of this element
\draw[black,dashed] (4.5,0.1) circle (0.8cm);

\draw[->] (5.1,0.25) -- ++(1.7,0.5) 
node [above,midway,sloped] {\scriptsize \,\,parameter};
\draw[->] (5.1,-0.25) -- ++(1.6,-0.5)
node [above,midway,sloped] {\scriptsize geometry};

%\crva(7,0.55);\crva(7,0.7);\crva(7,0.85);
\crva(7,0.3+0.75);
\crvaa(7,0.3+0.5);
\crvaaa(7,0.3+0.25);

\draw (6.8,-1.25) -- ++(1,0) -- ++(0,1) -- ++(-1,0) -- cycle;

\def\crs(#1,#2){
\draw (#1,#2) -- ++(0.7,0.5);
\draw (#1,#2+0.5) -- ++(0.7,-0.5);
\draw (#1,#2+0.5) -- ++(0.7,0.5);
\draw (#1,#2+1) -- ++(0.7,-0.5);
\node at (#1+1,#2) {...};
\node at (#1+1,#2+0.5) {...};
\node at (#1+1,#2+1) {...};
\draw (#1+1.3,#2) -- ++(0.7,0.5);
\draw (#1+1.3,#2+0.5) -- ++(0.7,-0.5);
\draw (#1+1.3,#2+0.5) -- ++(0.7,0.5);
\draw (#1+1.3,#2+1) -- ++(0.7,-0.5);
}
\def\nnlay(#1,#2){
\draw[black] (#1,#2) circle (4pt);
\draw[black] (#1,#2+0.5) circle (4pt);
\draw[black] (#1,#2+1) circle (4pt);
\draw[black] (#1,#2+1.5) circle (4pt);
\draw[black] (#1,#2+2) circle (4pt);
\draw[black] (#1,#2+2.5) circle (4pt);
\draw[black] (#1,#2+3) circle (4pt);
\draw[black] (#1,#2+3.5) circle (4pt);
\draw[black] (#1,#2+4) circle (4pt);

\crs(#1,#2)
\crs(#1,#2+1)
\crs(#1,#2+2)
\crs(#1,#2+3)
\node at (#1+1,#2-0.5) {\scriptsize ANN};

\draw[black] (#1+2,#2) circle (4pt);
\draw[black] (#1+2,#2+0.5) circle (4pt);
\draw[black] (#1+2,#2+1) circle (4pt);
\draw[black] (#1+2,#2+1.5) circle (4pt);
\draw[black] (#1+2,#2+2) circle (4pt);
\draw[black] (#1+2,#2+2.5) circle (4pt);
\draw[black] (#1+2,#2+3) circle (4pt);
\draw[black] (#1+2,#2+3.5) circle (4pt);
\draw[black] (#1+2,#2+4) circle (4pt);
}

\draw (8.2,1) -- ++(1.1,1) -- ++(0,-4) -- ++(-1.1,1) -- cycle;
\node at (8.75,0) {\scriptsize encoder};

\nnlay(9.5,-2)

\draw (11.7,2) -- ++(1.1,0) -- ++(0,-4) -- ++(-1.1,0) -- cycle;
\node at (12.25,0) {\scriptsize decoder};

\inout(13.25,-1.5,1.3)
\node at (13.85,-1.9) {\scriptsize in2out};

\insol(13.25,0.6,1.3)
\node at (13.85,0.4) {\scriptsize in2sol};

\end{tikzpicture}
\caption{Element learning - architecture. A finite element mesh is broken down element by element. For each element, a neural network takes the element geometry and PDE parameters as inputs, and returns the in2sol and in2out operators, \eqref{eq:in_out} and \eqref{eq:in_sol}, as outputs. Then, beyond this schematic diagram, the in2out operators of neighboring elements are systematically coupled together to assemble a global solver (see Figure~\ref{fig:hdg_diag}).}
\label{fig:elem_lrn}
\end{figure}

Certain benefits arise from defining the operators in \eqref{eq:in_out} and \eqref{eq:in_sol} and the nonlinear mapping in \eqref{eq:nn} on an element, rather than over the entire domain. For instance, in working on an individual element rather than the entire domain,
%On each element, since 
there is a significant reduction in the variations from the parameter $\sigma$, and in the complexity of geometry. Hence it is reasonable to anticipate that the neural network can provide a more reliable and accurate approximation to this non-linear mapping \eqref{eq:nn} on an element, compared to a neural-network-based approximation of this mapping on the whole domain. 

Element learning can be regarded as a framework that bridges the data-driven machine learning approach and the traditional finite element approach. To see this more clearly, note that in the extreme cases, such as when we have only a single element covering the whole domain, the methodology simplifies into a purely data-driven approach where the data are generated by a finite element method with a single element (with potentially high-order approximations leading it to a spectral-like method). Conversely, if we minimize the mesh size or decrease the polynomial degree, causing a significant reduction of the variations of the discretization of the coefficients distribution, the methodology approximates a traditional finite element method.

%We would like to highlight 
Some advantageous features of element learning 
%compared to purely data-driven machine learning approaches 
are the following:

\begin{itemize}
    \item The methodology can inherit numerous beneficial properties from classical finite element methods, such as proficient handling of complex geometry and boundary conditions, as it operates within the traditional finite element framework.
    \item A computational speedup is expected for element learning in comparison to traditional finite element-type methods, since element learning brings a memory reduction via hybridization (as in, e.g., HDG) and a faster creation of the in2out and in2sol operators via machine learning.
    \item The methodology doesn't require retraining when geometry changes or mesh size is refined or coefficients are modified. Training is only necessary for a single element with various geometrical and coefficient variations. Subsequently, this element can be used freely to solve any domain geometry and coefficient distribution.   
    \item Compared to many data-driven approaches, element learning is expected to be more reliable. This is due to the considerable reduction in geometric complexity and coefficient variations within each element. As the mesh is further refined, the reliability of the element learning approach increases.
    \item The methodology facilitates a natural coupling between data-driven and classical finite element approaches. 
    For instance, one could use the element learning solution as a preconditioner/initial condition to accelerate the convergence of a classical finite element method.
    \item Element learning has a good potential to be used for massive-scale parallel computing on GPUs since the forward propagation of the neural network for the in2sol and in2out operators on each element are completely independent processes.
\end{itemize}

% A significant challenge faced when numerically solving this equation stems from its high dimensionality. With three spatial dimensions, two angular dimensions, and one dimension for light frequency, conventional methods such as finite difference, finite volume, or finite element, may struggle to provide adequately precise solutions within the bounds of limited computational resources or stringent time constraints.

%\subsection{Organization of the paper}
The organization of the remainder of the paper is as follows.
In Section \ref{sec:num_method}, we introduce the element learning method for radiative transfer equations. 
Then in Section \ref{sec:num_exp}, we present numerical experiments to test the performance of the element learning methods in various scenarios, by comparing it with classical finite element-type methods such as DG and HDG methods.
A concluding discussion is in Section \ref{sec:conclusion}.

\section{Element learning for radiative transfer}
\label{sec:num_method}

In this section, we begin with subsection \ref{sec:rt_eqns}, in which we introduce the radiative transfer or radiation transport equation that is considered in this paper.

Then, in the subsection \ref{sec:HDG}, we show how to obtain a discretization of the in2out and in2sol operators by hybridizable discontinuous Galerkin method, which lays the foundation for element learning. This subsection include three parts: discontinuous Galerkin (DG), hybridizable discontinuous Galerkin (HDG), and the calculation for the in2out and in2sol operators. 

Finally, in the subsection \ref{sec:nn} we introduce the element learning accelerating method, including the neural network structure and the data generation algorithm.

\subsection{Radiative transfer}
\label{sec:rt_eqns}
Let $\Omega$ be a spatial domain and $S$ be the unit sphere in $\mathbb R^d$. We first introduce the inflow and outflow boundaries
\begin{align*}
    \Gamma^- :=\{(\mb x,\mb s)\in\pp\Omega\times S:\, \mb n(\mb x)\cdot\mb s\le 0\},\quad
    \Gamma^+ :=\{(\mb x,\mb s)\in\pp\Omega\times S:\, \mb n(\mb x)\cdot\mb s\ge 0\},
\end{align*}
where $\mb n(\mb x)$ is the outward-pointing normal vector at $\mb x\in\pp\Omega$. We consider the following radiative transfer equation
\begin{subequations}
\label{eq:RT_PDE}
\begin{alignat}{5}
    \mb s\cdot\nabla u + \sigma_e u - \sigma_s \int_S p(\mb s,\mb s')u(\mb s')ds' &= f &\qquad&\mr{in}\ \Omega,\\
    u &= g && \mr{on}\ \Gamma^-,
\end{alignat}
\end{subequations}
where the solution $u(\mb x,\mb s)$ represents the radiant intensity at $\mb x$ in the direction $\mb s\in S$, $\sigma_e(\mb x)$ and $\sigma_s(\mb x)$ are the extinction and scattering coefficients, respectively, $f(\mb x,\mb s)$ and $g(\mb x,\mb s)$ are the source term and the inflow radiation, respectively. Here $p(\mb s,\mb s')$ is the scattering phase function, for which we consider the Henyey-Greenstein phase function
\begin{align*}
    p(\mb s,\mb s')=\frac{1-g_\mr{asym}^2}{c(1+g_\mr{asym}^2-2g_\mr{asym}\cos\mr{ang}(\mb s,\mb s'))^{3/2}},
\end{align*}
where $\mr{ang}(\mb s,\mb s')$ denotes the angle between $\mb s$ and $\mb s'$, $g_\mr{asym}$ is the asymmetric parameter taking values in $[0,1]$, and the constant $c$ is chosen such that $\int_Sp(\mb s,\mb s')ds'=1$. For isotropic scattering we have $g=0$, while for the scattering of short-wave (solar) radiation in water clouds, $g$ can take values from 0.8 to 0.9 \cite{slingo1989gcm}. 

\subsection{Discretizations of the in2out and int2sol operators}
In this subsection, we show how to obtain a discretization for the in2out and the in2sol operators by hybridizable discontinuous Galerkin (HDG) methods. Since this is likely the first paper that introduces HDG methods for radiative transfer, we shall introduce this method with some details in the following three parts: Discontinuous Galerkin (DG) methods, HDG methods, and the implementation of HDG for the discretization of the in2out and in2sol operators.

\label{sec:HDG}
\subsubsection{Discontinuous Galerkin (DG) methods}
We begin by partitioning the spatial domain $\Omega$ into a collection of polyhedral elements $\mc T_h$, and partitioning the unit sphere $S$ into a collection of angular elements $\mc T_h^a$. 
We shall assume the partitions are conforming and shape-regular; these are standard assumptions for finite element-type methods and we refer to \cite{BrSc:2008} for more details.

Let $\mc E_h$ be the collection of all faces of $\mc T_h$, and let $\mc E_h^0$ and $\mc E_h^\Gamma$ be the collection of the interior faces ($F\nsubset\pp\Omega$) and boundary faces ($F\subset\pp\Omega$), respectively. 
We also define the skeleton of the mesh as $\mc E:=\cup_{F\in\mc E_h} F$.
For each element $K\in\mc T_h$, we denote by $\mc E_K$  the collection of the faces of $K$. We require that the angular partition $\mc T_h^a$ respects the spatial partition $\mc T_h$ in the sense that for any $K\times K^a\in\mc T_h\times\mc T_h^a$ and a given face $F\in\mc E_K$, the angular element $K^a$ is either entirely inflow or entirely outflow with respect to the face $F$ and the element $K$.

We denote by $\mc P_k(K)$  the polynomial space with degree $k$ supported on the element $K$. Now we can introduce the following DG approximation spaces:
\begin{align*}
    V_h:=\{u_h\in L^2(\Omega\times S):\, u_h\big|_{K\times K^a}\in\mc P_k(K)\otimes\mc P_{k^a}(K^a)\quad\forall\,K\times K^a\in\mc T_h\times\mc T_h^a\}.
\end{align*}
Namely, $u_h$ is piece-wise polynomial on each spatial-angular element $K\times K^a\in\mc T_h\times\mc T_h^a$. 

For the DG methods, we seek for $u_h\in V_h$ such that
\begin{subequations}
\begin{align}
    \nonumber
    \sum_{K\in\mc T_h}\sum_{K^a\in\mc T_h^a}\Big(\int_{K^a}\int_{\pp K} (\mb s\cdot\mb n)\widehat{u}_h v - \int_{K^a}\int_K u_h\mb s\cdot\nabla v + \int_{K^a}\int_K\sigma_e u_h v&\\
    \label{eq:dg}
    -\int_{K^a}\int_K\sigma_s(\mb x)\int_Sp(\mb s,\mb s')u_h(\mb x,\mb s')ds'v(\mb x,\mb s)dxds\Big) &= \sum_{K\in\mc T_h}\sum_{K^a\in\mc T_h^a}\int_{K^a}\int_K f v,
\end{align}
hold for all $v\in V_h$,
where the numerical flux quantity $\widehat u_h$ is defined as follows:
\begin{align}\label{eq:dg_flux}
    \widehat{u}_h(\mb x,\mb s) = \left\{
    \begin{array}{ll}
        u_h(\mb x,\mb s) &  \text{if}\quad \mb n_F\cdot\mb s\ge0,\\
        u_h^\mr{nbr}(\mb x,\mb s) & \text{if}\quad \mb n_F\cdot\mb s\le 0\text{ and }F\in\mc E_h^0,\\
        g(\mb x,\mb s) & \text{if}\quad \mb n_F\cdot\mb s\le 0\text{ and }F\in\mc E_h^\Gamma,
    \end{array}
    \right.
\end{align}    
\end{subequations}
where $u_h^\mr{nbr}$ is the value of $u_h$ from the neighbouring elements.
Namely, we choose the upwind flux for $\widehat u_h$. 
We refer to \cite{du2023fast,du2023inverse} and references therein for more details on the DG discretization of the radiative transfer equation.

\subsubsection{Hybridization}
For hybridizable DG (HDG) methods, we introduce another variable $\widehat u_h$ living on the skeleton $\mc E$. This variable (also known as hybrid unknown) can be regarded as a Lagrange multiplier of enforcing the consistency condition. To begin with, let us introduce its corresponding approximation space:
\begin{align*}
    \widehat V_h:=\{\widehat u_h\in L^2((\cup_{F\in\mc E_h} F)\times S):\, \widehat u_h\big|_{F\times K^a}\in\mc P_k(F)\otimes\mc P_{k^a}(K^a)\quad\forall\,F\times K^a\in\mc E_h\times\mc T_h^a\}.
\end{align*}
Namely, $\widehat u_h$ is a piece-wise polynomial on the skeleton $\mc E$, where each piece of the skeleton is a face $F\in\mc E_h$.
We shall also introduce
\begin{align}
    \widehat V_h^g:=\{\widehat u_h\in\widehat V_h:\, \widehat u_h = \mr P_M g\quad\forall F\times K^a\in\mc E_h\times\mc T_h^a\quad \text{s.t.}\quad
    F\subset\pp\Omega,\, \mb n_F\cdot K^a\le 0\},
\end{align}
where $\mr P_M$ is the orthogonal projection from the space $L^2(F\times K^a)$ to the space $\mc P_k(F)\otimes P_{k^a}(K^a)$. This space can be regarded as projecting the boundary data $g$ to the hybrid unknown space $\widehat V_h$.

For HDG methods, we seek $(u_h,\widehat u_h)\in V_h\times \widehat V_h^g$ such that
\begin{subequations}
\label{eq:hdg_loc_glb}
\begin{align}
    \nonumber
    \sum_{F\in\mc E_K}\sum_{K^a\in\mc T_h^a\cap(K^a\cdot\mb n_F\ge0)}\int_{K^a}\int_F (\mb s\cdot\mb n_F)u_h v&\\
    \nonumber
    +\sum_{F\in\mc E_K}\sum_{K^a\in\mc T_h^a\cap(K^a\cdot\mb n_F\le0)}\int_{K^a}\int_F (\mb s\cdot\mb n_F)\widehat{u}_h v&\\
    \nonumber
    +\sum_{K^a\in\mc T_h^a}\left( - \int_{K^a}\int_K u_h\mb s\cdot\nabla v + \int_{K^a}\int_K\sigma_e u_h v\right)&\\
    \label{eq:hdg_local}
    -\sum_{K^a\in\mc T_h^a}\int_{K^a}\int_K\sigma_s(\mb x)\int_Sp(\mb s,\mb s')u_h(\mb x,\mb s')ds'v(\mb x,\mb s)dxds &= \sum_{K^a\in\mc T_h^a}\int_{K^a}\int_K f v\quad\forall v\in V_h,
\end{align}    
holds for all $K\in\mc T_h$. Note that \eqref{eq:hdg_local} can be solved element-wise for each $K\in\mc T_h$, if given $\widehat u_h$ and $f$ as the input data. Namely, by solving \eqref{eq:hdg_local}, we can element-wise express the solution $u_h$ as a function of $\widehat u_h$ and $f$ by $u_h=u_h(\widehat u_h,f)$. 
The global solver is nothing but enforcing the following consistency condition:
\begin{align}\label{eq:hdg_global}
    \sum_{K\in\mc T_h}\sum_{F\in\mc E_K}\sum_{K^a\in\mc T_h^a\cap(K^a\cdot\mb n_F\ge 0)}\int_F\int_{K^a} \left(\widehat{u}_h-u_h(\widehat u_h,f)\right)\eta = 0\quad\forall\eta\in\widehat V_h^0.
\end{align}
\end{subequations}

\subsubsection{The in2out and in2sol operators}
\label{sec:implem_hdg}
In this part we show how to derive a discretization for the in2out and in2sol operators by the HDG formulation \eqref{eq:hdg_loc_glb}.

Let $\varphi_j^{K\times K^a}$ be the basis of $\mc P_k(K)\otimes P_{k^a}(K^a)$, and $\psi_j^{F\times K^a}$ be the basis of $\mc P_k(F)\otimes\mc P_{K^a}(K^a)$. 
Then, the solution $u_h$ and the hybrid unknown $\widehat u_h$ can be represented as the following linear combinations:
\begin{align*}
    u_h =\sum_{K\in\mc T_h}\sum_{K^a\in\mc T_h^{a}} \sum_{j}u_{j}^{K\times K^a}\, \varphi_{j}^{K\times K^a},\qquad
    \widehat u_h = \sum_{F\in\mc E_h}\sum_{K^a\in\mc T_h^{a}} \sum_{j}\widehat u_{j}^{F\times K^a}\, \psi_{j}^{F\times K^a}.    
\end{align*}
We can rewrite \eqref{eq:hdg_local} into a matrix form as follows:
\begin{align}\label{eq:hdg_local_mat}
    B^K[u_h] + \widehat B^K[\widehat u_h] - C^K[u_h] + M^K[u_h]
    - S^K[u_h] = [f],
\end{align}
where $[u_h]$, $[\widehat u_h]$ and $[f]$ represents the DOFs of $u_h$, $\widehat u_h$ and $f$, respectively. Namely,
\begin{align*}
    [u_h]_{K,K^a,i} = u_i^{K\times K^a},\quad
    [\widehat u_h]_{F,K^a,i} = \widehat u_i^{F\times K^a},\quad
    [f]_{K,K^a,i} = \int_K\int_{K^a}f\varphi_i^{K\times K^a}.
\end{align*}
The matrix terms of \eqref{eq:hdg_local_mat} are assembled as follows (let the test function $v = \varphi_{i}^{K\times K^a}$):
\begin{subequations}
\label{eq:hdg_local_all_sub_mat}
\begin{align}
    (B^K[u_h])_{K,K^a,i} &=\sum_{j}u_{j}^{K\times K^a}
    \left(
    \sum_{F\in\mc E_K}
    \mathbbm 1_{K^a\cdot\mb n_F\ge0}\int_{K^a}\int_F (\mb s\cdot\mb n_F)\, \varphi_{j}^{K\times K^a} \varphi_{i}^{K\times K^a}
    \right),\\
    (\widehat B^K[\widehat u_h])_{K,K^a,i}&=
    \sum_{F\in\mc E_K}\sum_{j}\widehat u_{j}^{F\times K^a}\left(\mathbbm 1_{K^a\cdot\mb n_F\le0}\int_{K^a}\int_F (\mb s\cdot\mb n_F)  \psi_{j}^{F\times K^a}\varphi_i^{K\times K^a}\right),\\
    (C^K[u_h])_{K,K^a,i}&=\sum_{j}u_{j}^{K\times K^a}\left(\int_{K^a}   \int_K  \varphi_{j}^{K\times K^a}\mb s\cdot\nabla \varphi_i^{K\times K^a}\right),\\
    (M^K[u_h])_{K,K^a,i}&=\sum_ju_j^{K\times K^a}\left(\int_{K^a}\int_K\sigma_e\varphi_j^{K\times K^a}\varphi_i^{K\times K^a}\right),\\
    (S^K[u_h])_{K,K^a,i}&=\sum_{K_{*}^a\in\mc T_h^a}\sum_{j}u_{j}^{K\times K_{*}^a}
    \left(\int_K\sigma_s(\mb x)\int_{K^a}\int_{K_*^a}p(\mb s,\mb s')
    \varphi_{j}^{K\times K_{*}^a}(\mb x,\mb s')
    \varphi_i^{K\times K^a}(\mb x,\mb s)ds'dsdx\right).
\end{align}    
\end{subequations}
We leave the details on how the above formulations are derived in the Appendix \ref{sec:app_hdg_mat}.

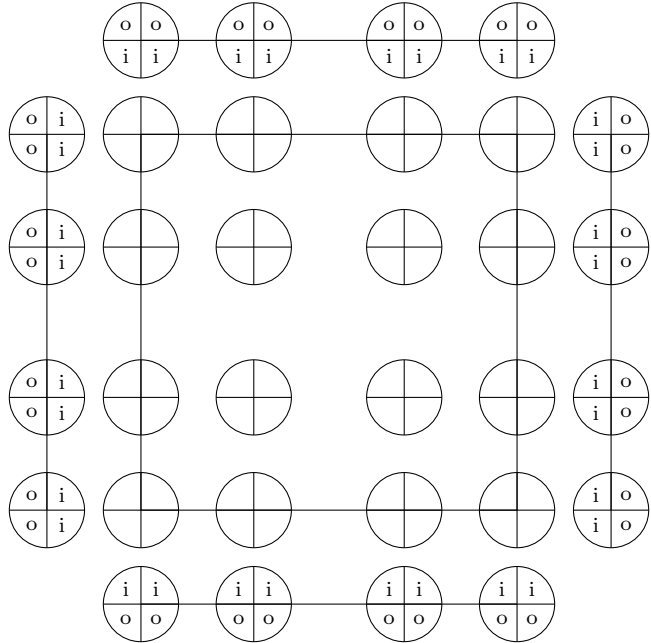
\begin{figure}[H]
\centering
\begin{tikzpicture}
\draw (0,0) -- ++(5,0) -- ++(0,5) -- ++(-5,0) -- cycle;

\newcommand{\crccrs}[2]{
\draw (#1,#2) circle (0.5);
\draw (#1-0.5,#2) -- (#1+0.5,#2);
\draw (#1,#2-0.5) -- (#1,#2+0.5);
}

\newcommand{\crccrsLin}[2]{
\crccrs{#1}{#2}
\crccrs{#1+1.5}{#2}
\crccrs{#1+3.5}{#2}
\crccrs{#1+5}{#2}
}

\crccrsLin{0}{0}
\crccrsLin{0}{1.5}
\crccrsLin{0}{3.5}
\crccrsLin{0}{5}

\def\ooii(#1,#2){
\node at (#1-0.2,#2+0.2) {\scriptsize o};
\node at (#1+0.2,#2+0.2) {\scriptsize o};
\node at (#1-0.2,#2-0.2) {\scriptsize i};
\node at (#1+0.2,#2-0.2) {\scriptsize i};
}
\def\iioo(#1,#2){
\node at (#1-0.2,#2+0.2) {\scriptsize i};
\node at (#1+0.2,#2+0.2) {\scriptsize i};
\node at (#1-0.2,#2-0.2) {\scriptsize o};
\node at (#1+0.2,#2-0.2) {\scriptsize o};
}
\def\oioi(#1,#2){
\node at (#1-0.2,#2+0.2) {\scriptsize o};
\node at (#1+0.2,#2+0.2) {\scriptsize i};
\node at (#1-0.2,#2-0.2) {\scriptsize o};
\node at (#1+0.2,#2-0.2) {\scriptsize i};
}
\def\ioio(#1,#2){
\node at (#1-0.2,#2+0.2) {\scriptsize i};
\node at (#1+0.2,#2+0.2) {\scriptsize o};
\node at (#1-0.2,#2-0.2) {\scriptsize i};
\node at (#1+0.2,#2-0.2) {\scriptsize o};
}

\def\crccrsLintp(#1,#2){
\crccrs{#1}{#2}
\crccrs{#1+1.5}{#2}
\crccrs{#1+3.5}{#2}
\crccrs{#1+5}{#2}
\draw (#1,#2) -- ++(5,0);
\ooii(#1,#2)
\ooii(#1+1.5,#2)
\ooii(#1+3.5,#2)
\ooii(#1+5,#2)
}
\def\crccrsLinbt(#1,#2){
\crccrs{#1}{#2}
\crccrs{#1+1.5}{#2}
\crccrs{#1+3.5}{#2}
\crccrs{#1+5}{#2}
\draw (#1,#2) -- ++(5,0);
\iioo(#1,#2)
\iioo(#1+1.5,#2)
\iioo(#1+3.5,#2)
\iioo(#1+5,#2)
}
\def\crccrsLinlf(#1,#2){
\crccrs{#1}{#2}
\crccrs{#1}{#2+1.5}
\crccrs{#1}{#2+3.5}
\crccrs{#1}{#2+5}
\draw (#1,#2) -- ++(0,5);
\oioi(#1,#2)
\oioi(#1,#2+1.5)
\oioi(#1,#2+3.5)
\oioi(#1,#2+5)
}
\def\crccrsLinrt(#1,#2){
\crccrs{#1}{#2}
\crccrs{#1}{#2+1.5}
\crccrs{#1}{#2+3.5}
\crccrs{#1}{#2+5}
\draw (#1,#2) -- ++(0,5);
\ioio(#1,#2)
\ioio(#1,#2+1.5)
\ioio(#1,#2+3.5)
\ioio(#1,#2+5)
}

\crccrsLintp(0,6.25)
\crccrsLinbt(0,-1.25)
\crccrsLinlf(-1.25,0)
\crccrsLinrt(6.25,0)
\end{tikzpicture}
\caption{A visual representation of the DOFs of the solution $u_h$ and $\widehat u_h$. A single $Q_3$ spectral element ($p_x=p_y=3$) is shown with a circle at each of the element's $4\times 4=16$ spatial quadrature points. Each circle also represents the angular coordinate, as $4$ uniformly-partitioned $P_0$ angular elements ($N_a=4$ and $p_a=0$) at each spatial quadrature point. On the skeleton along the boundary of the spatial element, the letters `i' and `o' represent the inflow and the outflow DOFs of $\widehat u_h$, respectively.}
\label{fig:hdg_dofs_visual}
\end{figure}

By inverting \eqref{eq:hdg_local_mat}, we obtain
\begin{align}\label{eq:hdg_local_mat_sol}
    [u_h] = A_{i2u}^K[\widehat u_h] + f_u^K,
\end{align}
where
\begin{subequations}
\label{eq:hdg_local_f2u_i2u}
\begin{align}
\label{eq:hdg_local_i2u}
    A_{i2u}^K &:= -(B^K-C^K+M^K-S^K)^{-1}\widehat B^K,\\
\label{eq:hdg_local_f2u}
    f_u^K &:= (B^K-C^K+M^K-S^K)^{-1}[f].
\end{align}    
\end{subequations}
In the above formulation, the matrix $A_{i2u}^K$ maps the inflow radiation $\widehat u_h$ to the corresponding part of the numerical solution $u_h$, this represents a discretization of the {\sl in2sol} operator, and
$f_u^K$ is the part of the solution $u_h$ that corresponds to the forcing $f$. Based on them, we immediately obtain
\begin{subequations}
\label{eq:hdg_local_f2o_i2o}
\begin{alignat}{5}
\label{eq:hdg_local_i2o}
    A_{i2o}^K &= R_K^oA_{i2u}^K && \text{(inflow to outflow radiation)},\\
    \widehat f_{u}^K &= R_K^of_{u}^K &\quad& \text{(forcing on skeleton)},
\end{alignat}    
\end{subequations}
where $R_K^o$ denotes a restriction of the DOFs of $u_h$ (on the element $K$) to the outflow DOFs of $\widehat u_h$ on $K$; see Figure \ref{fig:hdg_dofs_visual} for a visualization of the DOFs on a cubic element. Here $A_{i2o}^K$ represents a discretization of the {\sl in2out} operator.

By \eqref{eq:hdg_global} and \eqref{eq:hdg_local_f2o_i2o}, we have
\begin{align}\label{eq:hdg_global_mat}
    \sum_{K\in\mc T_h} \left(
    \widehat R_K^{o}[\widehat u_h]-
    A_{i2o}^K[\widehat u_h]\right)=\sum_{K\in\mc T_h}\widehat f_{u}^K.
\end{align}
Here, the matrix $\widehat R_K^o$ restricts the DOFs of $\widehat u_h$ to the outflow DOFs of $\widehat u_h$ on element $K$. Note that the above system \eqref{eq:hdg_global_mat} only involves the hybrid unknown $\widehat u_h$.

Let us conclude this subsection by the following Algorithm \ref{algrm:hdg}.
\begin{algorithm}[H]
\caption{HDG implementation}
\label{algrm:hdg}
\begin{algorithmic}[1]
\State {\bf Creation of local solver}: for each element $K\in\mc T_h$, solve the local system \eqref{eq:hdg_local_f2u_i2u} to obtain $A_{i2u}^K$ (discretization of the Green's function) and $f_u^K$.
\State Assemble the global system \eqref{eq:hdg_global_mat} using $A_{i2o}^K$ and $\widehat f_u^K$.
\State {\bf Global solve}: solve the global system \eqref{eq:hdg_global_mat} to obtain $\widehat u_h$.
\State For each element $K\in\mc T_h$, recover the solution $u_h$ by using \eqref{eq:hdg_local_mat_sol}, which needs the hybrid solution $\widehat u_h$ obtained in Step 3, and the matrix $A_{i2u}^K$ and the term $f_u^K$ obtained in Step 1.
\end{algorithmic}
\end{algorithm}

\subsection{Solvers for the linear algebraic systems}
\label{sec:solver}
{Many solvers are in use for radiative transfer, such as source iteration, sweeps, diffusion synthetic acceleration, preconditioned GMRES, multigrid, or combinations among them; we refer to \cite{adams2002fast,azmy2010advances} for a review on this topic and \cite{morel1991angular,warsa2004krylov,chacon2017multiscale,dolz2022robustly} for more. 
%Different solvers are used in different applications, including atmospheric science, neutron transport, medical imaging, astrophysics, and others. 
Performance of the solver can vary depending on the problem itself (e.g., strong or weak scattering), the numerical discretizations (e.g., discrete ordinate, spherical harmonic, finite volume, discontinuous Galerkin), serial or parallel implementations, software being used (e.g., Matlab, C, Python, Fortran), and
the computer architecture (CPU or GPU, parallelization support, memory speed).

Here, for the DG method, we invert the system \eqref{eq:dg} using a preconditioned GMRES method, similar to
\cite{patton2002application,ReBaHi:2006}, where the preconditioner is obtained by inverting the advection-extinction subsystem with UMFPACK \cite{davis2004algorithm}.
We have examined several test cases
and observed in \cite{du2023fast} that this solver 
{scales near optimally as $\mc O(N^{1.1})$ where $N$ is the total number of degrees of freedom, and it}
outperforms solver choices such as source iteration or directly applying UMFPACK to the whole system. 
For the HDG method, note that we have two systems to invert, namely the local solver and the global solver; see Algorithm \ref{algrm:hdg}. This differs from the DG method, for which we only need to invert the global system \eqref{eq:dg}. The local systems of HDG are inverted by UMFPACK. 
Then, the global system is inverted by a GMRES method \cite{saad1986gmres}. 
We remark that GMRES is a Krylov subspace method and supports matrix-free implementation, which would be an interesting future direction.
%For HDG-EL, the global system is also inverted by GMRES for a fair comparison.
}

% {It should be noted that our selection of solvers is not claimed to be the most optimized. However, they are extensively used across various applications. For instance, UMFPACK is recognized as a general-purpose solver (and also the default solver in Matlab) for unsymmetric sparse linear systems \cite{davis2004algorithm}. Also, GMRES is a widely adopted Krylov subspace iterative solver for unsymmetric systems \cite{saad1986gmres}, known for its compatibility with matrix-free implementations \cite{johan1991globally}.

% Moreover, it is important to point out that achieving a fully optimized solver typically requires a complex optimization process that spans the PDE itself, the corresponding numerical algorithms, software, and hardware. This complexity often results in a lack of generality, presenting coding challenges and/or necessitating significant domain-expert efforts for different problems.
% In this context, one of our objectives in introducing element learning is to offer a systematic and universal strategy for enhancing the efficiency of these general solvers, irrespective of the problem. This approach aims not only to accelerate computational processes but also to mitigate the extensive optimization efforts traditionally required, making solver application more accessible and less labor-intensive for researchers and practitioners alike.
% }

\subsection{Neural network}
\label{sec:nn}
The HDG method introduced in the previous subsection allows a reduction of the DOFs of $u_h$, which is defined on the whole DG space $V_h$, to the hybrid unknown $\widehat u_h$, which is  defined only on the skeleton space $\widehat V_h$. In addition, the higher order of the approximation we use, the more reduction that can be achieved, and so one can expect a faster convergence. 

However, a significant difficulty that hinders the use of very high order HDG methods is that the local system becomes much more expensive to invert as the polynomial degree increases. To see this, note that the local solver \eqref{eq:hdg_local_i2u} requires one to solve for {\sl all} possible combinations of the inflow boundary conditions. As a result, as the polynomial degree increases, the dimension of the inflow radiation space (the number of columns of $\widehat B^K$ in \eqref{eq:hdg_local_i2u}) increases dramatically. Considering that the system size (the number of rows of $\widehat B^K$ in \eqref{eq:hdg_local_i2u}) also increases, it very quickly becomes difficult and potentially infeasible to apply HDG methods as the polynomial degree becomes large. 
We also refer to \cite{woopen2014comparison}, where it is shown that a majority of the computational cost comes from inverting the local systems to assemble the global matrix.

To mitigate this issue, in this section, we aim to leverage a neural network to accelerate the local system inversion \eqref{eq:hdg_local_f2u_i2u}. 
Note that \eqref{eq:hdg_local_f2u} requires the inversion for only one instance $[f]$, so a major computational cost comes from \eqref{eq:hdg_local_i2u}, which requires the inversion for multiple instances where each instance is a column vector of $\widehat B^K$. 
So here we focus on using neural network to accelerate the calculation for $ A_{i2u}^K$ in \eqref{eq:hdg_local_i2u}. 

We remark that $A_{i2u}^K$ depends on the parameter $(\sigma,K)$, and what we aim to approximate is not $A_{i2u}^K$ for some particular $(\sigma,K)$. Instead, we aim to recover the full mapping $S_{DG}:(\sigma,K)\rightarrow A_{i2u}^K$, which takes $(\sigma,K)$ as input, and returns the corresponding Green's function discretization $A_{i2u}^K$ as output. Since $S_{DG}$ is non-linear, we cannot use linear regression. This justifies the use of neural networks with non-linear activation functions. Numerical evidences will be presented in Section \ref{sec:num_elem_train}.

One way is to use a neural network to directly approximate the operator $S_\mr{DG}:(\sigma,K) \rightarrow A_{i2u}^K$ with $A_{i2u}^K$ as its output. However, in many cases, the full solution $u_h$ from $A_{i2u}^K$ is not needed, and only some partial statistics of the solution $u_h$ (e.g. mean intensity, radiative heating rate) are needed. For instance, here we consider the case of recovering the mean intensity: 
\begin{align}
\label{eq:ang_ave}
    u_h^m(\mb x) = \frac{1}{|S|}\int_S u_h(\mb x,\mb s)ds.
\end{align}
In this case, we don't need to calculate the full matrix $A_{i2u}^K$ but only a small fraction of the statistics of this matrix. 
To be more specific, we only need the following two matrices: (1) $A_{i2o}^K$ (defined in \eqref{eq:hdg_local_i2o}) for the assembling of the global system \eqref{eq:hdg_global_mat}, and (2) $A_{i2m}^K$ for recovering the mean intensity of $u_h$. Here $A_{i2m}^K:=R_{u2m}^KA_{i2u}^K$ and $R_{u2m}^K$ is the matrix form of the above angular-averaging operator defined in \eqref{eq:ang_ave}. Considering this, here we aim to construct a neural network (NN) to approximate the following operator
\begin{align}\label{eq:param2sol}
    \widetilde S_\mr{DG}: (\sigma,K)\rightarrow (A_{i2o}^K,A_{i2m}^K),
\end{align}
where $A_{i2o}$ is for inter-element communication and for assembling the global system \eqref{eq:hdg_global_mat} and $A_{i2m}$ is for local solution recovery for the mean intensity $u_h^m$.

{
Note that with this NN, the local solver creation process in Algorithm \ref{algrm:hdg} can be replaced by machine learning approaches for acceleration. As a result, for element learning, the only system we need to invert is the global system which shares the same structure as the global system of a standard HDG method. 
}

In the next subsection, we explain how to construct neural networks to the approximation of $\widetilde S_\mr{DG}$ in more details.

\subsubsection{Networks architecture}
\label{sec:net_struct}
To use a NN to approximate the operator $\widetilde S_\mr{DG}$ in \eqref{eq:param2sol}, we begin by considering the parameterization of the optical coefficient $\sigma$ and the element geometry $K$.
For the conciseness of the presentation and the convenience of the implementation, we consider the case when $d=2$ and $K$ is a rectangular element. 
The generalization to $d=3$ and other types of element (e.g. triangular elements) follow the same ideas.
We consider the case when the angular space $S$ is a unit circle, so the radiation direction $\mb s$ can be characterized by a single parameter, namely the angle $\theta\in[0,2\pi]$, by $\mb s=(\cos\theta,\sin\theta)$.
Also, we consider a fixed angular discretization of $[0,2\pi]$ by a uniform partition $\mc T_h^a$ with $N_a$ elements and polynomial degree $p_a$.

To consider the parameterization of element geometry $K$,
let $\widehat K = [-1,1]^2$ be the reference element, in which case the geometry of $K$ is determined by the following push-forward map:
\begin{align*}
    F_K: \widehat K &\rightarrow K,\\
    (\hat x,\hat y)&\mapsto (a_{11}^K\hat x+a_{12}^K\hat y,a_{21}^K\hat x+a_{22}^K\hat y)+(b_1^K,b_2^K).
\end{align*}
Namely, the geometry of $K$ is parameterized by the parameters $a_{ij}^K$ and $b_i^K$ with $i,j=1,2$. 
One could then proceed by including parameters $a_{ij}^K$ and $b_i^K$ into the input layer of the NN. However,
in many cases, it is straightforward to derive the change of the solution under a push-forward map. 
{So in these cases, we don't need to insert $a_{ij}^K$ and $b_i^K$ to the input layer of the neural network.}
For instance, if $K$ is a square element, then the push forward map $F_K$ becomes
$F_K(\hat x,\hat y)=h(\hat x,\hat y)+(b_1^K,b_2^K)$, where $h$ represents the size of the element $K$. In this case, it is straightforward to derive that solving \eqref{eq:RT_PDE} on $K$ is equivalent to solving the same equation on $\widehat K$ with re-scaled coefficients $h\sigma_e/2$, $h\sigma_s/2$, and forcing term $hf/2$. Thus, it suffices to consider the parameterization on the reference element $\widehat K$. For the following discussion, we focus on this case for simplicity.

{\color{black}
We remark that, if one were to use a mesh containing multiple element types (e.g. triangles, quadrilateral), two approaches could be considered.
For the first approach, one could use multiple neural networks so that each network learns the operator \eqref{eq:param2sol} for each element type. For the second approach, one could use a neural network to learn the operator \eqref{eq:param2sol} on a general polyhedral element, in which case the simple-geometry element can be regarded as a special case of the complex-geometry element; see, for instance \cite{DuSa_poly:2021,DiTi:2018}.
}

To parameterize $\sigma$, we consider a nodal representation by spectral element:
\begin{align}\label{eq:sigma_se}
    \sigma(x,y)\approx\sigma_h(x,y)= \sum_{i=1}^{p_x+1}\sum_{j=1}^{p_y+1}\sigma_{ij}\,\varphi_i(x)\varphi_j(y),
\end{align}
where $\varphi_i$ is the Lagrange polynomial associated with the $i$-th Legendre-Gauss-Lobatto (LGL) quadrature point on $[-1,1]$. 
Note that we need to parameterize both the extinction $\sigma_e$ and the scattering $\sigma_s$. 
As a result, the input layer of the network is a vector 
\begin{align}
    \mb x_0:=([\sigma_e],[\sigma_s])\in \mathbb R^{N_\mr{in}},    
\end{align}
where $N_\mr{in}=2(p_x+1)(p_y+1)$, and $[\sigma_e]$ and $[\sigma_s]$ are the two row vectors with their entries as the nodal representation of $\sigma_e$ and $\sigma_s$, respectively. For applications when $\sigma_e$ and $\sigma_s$ are correlated by $\sigma_s=\widetilde\omega\sigma_e$ where $\widetilde\omega$ is a given constant (e.g., short-wave radiation by cloud), we only need to parameterize one of these two coefficients so the input size becomes $N_\mr{in}=(p_x+1)(p_y+1)$.

We next consider the output layer. Note the matrix $A_{i2o}^K$ has the size of 
$2(p_x+p_y+2)\frac{N_a}{2}(p_a+1)\,\,
\times\,\,
2(p_x+p_y+2)\frac{N_a}{2}(p_a+1)$, and the matrix $A_{i2m}^K$ has the size
$(p_x+1)(p_y+1)\,\,
\times\,\,
2(p_x+p_y+2  )\frac{N_a}{2}(p_a+1)$.
As a result,
the output layer is a vector
\begin{align}
    \mb x_{N_l}:=([A_{i2o}^K],[A_{i2m}^K])\in \mathbb R^{N_\mr{out}}
\end{align}
where
$N_\mr{out}=(p_x+p_y+2)N_a(p_a+1)\left(
(p_x+p_y+2)N_a(p_a+1)+(p_x+1)(p_y+1)
\right)$, and $[A_{i2o}^K]$ and $[A_{i2m}^K]$ are the two row vectors collecting the entries of the matrices $A_{i2o}^K$ and $A_{i2m}^K$, respectively.

Finally, we construct a fully connected neural network from $\mb x_0$ to $\mb x_{N_l}$. Namely,
\begin{align*}
    \mb x_{i} = f_{i}(A_{i}\mb x_{i-1}+\mb b_{i}),\qquad  i=1,2,...,N_l,
\end{align*}
where $A_i$ is a linear transformation, $\mb b_i$ is a vector sharing the same size of $\mb x_{i}$, and $f_i$ is a non-linear activation function for $i=1,2,...,N_l-1$ and $f_{N_l}(\mb x)=\mb x$. We can also write the neural network in a compact form as follows:
\begin{align*}
    \mb x_{N_l} = \mr{NN}_{\bm\theta}(\mb x_0),\quad\mr{where}\quad
    \bm\theta = (\bm\theta_1,...,\bm\theta_{N_l})=(A_1,\mb b_1,...,A_{N_l},\mb b_{N_l}).
\end{align*}

\subsubsection{Data generation}
\label{sec:data_gen}
Now we consider the data generation that will be used to train the neural network $\mr{NN}_{\bm\theta}$. 
{\color{black}We remark that only the data generation on the reference element $\widehat K$ is considered here, since the operator on the physical element $K$ \eqref{eq:param2sol} can be derived based on the operator on the reference element $\widehat K$; see the second paragraph of Section \ref{sec:net_struct} for more details.}
For simplicity, we shall consider the case with a fixed single-scattering albedo $\tilde\omega$, which relates extinction and scattering by $\sigma_s = \tilde\omega\sigma_e$. In this case, it suffices to consider the data generation for only one optical coefficient instead of two, and we will use simply $\sigma_h$ to represent this coefficient.

We shall generate $\sigma_h$ according to a probability distribution $\mu_\sigma$. Then, the corresponding matrices $A_{i2o}^{\widehat K}(\sigma_h)$ and  $A_{i2m}^{\widehat K}(\sigma_h)$ will be calculated by using $\sigma_h$ and solving the local system on the reference element ${\widehat K}$; see Section \ref{sec:implem_hdg} for more details. We next explain how we define the probability distribution $\mu_\sigma$.

In \eqref{eq:sigma_se}, we define $\sigma_{ij}$ to be the nodal representation of $\sigma_h$ by spectral element. We can also express $\sigma_h$ in terms of modal representation:
\begin{align*}
    \sigma_h(x,y) = \sum_{m=0}^{p_x} \sum_{n=0}^{p_y} \widetilde\sigma_{mn}L_m(x)L_n(y),
\end{align*}
where $L_m$ is the $m$-th Legendre polynomial on $[-1,1]$. Note that $L_m$ with a higher index $m$ represents higher frequency signals. Using this fact, we shall generate the probability distribution $\mu_\sigma$ as described in Algorithm \ref{algrm:opt_data_gen}.

\begin{algorithm}[H]
\caption{Optical coefficient data generation}
\label{algrm:opt_data_gen}
\begin{algorithmic}[1]
\State Draw samples of the modal representation $\widetilde\sigma_{mn}$ with low-frequency bias:
    \begin{align*}
        \widetilde\sigma_{mn} = \exp{\left(-c_\mr{sm}\left((\frac{m}{p_x})^2-(\frac{n}{p_y})^2\right)\right)}(X_{mn}-\frac{1}{2}),
    \end{align*}
    where $X_{mn}$ are i.i.d uniform distribution on $[0,1]$, and $c_\mr{sm}$ determines the decaying speed of the high-frequency signals in the distribution $\sigma_h$ so a larger $c_\mr{sm}$ is likely to generate a smoother $\sigma_h$. This procedure is guided by the consideration that low-frequency bias exists universally in the physical world.
    
    \State Perform transformation from modal to nodal representations:
    \begin{align*}
        \sigma_{ij} = T_{\mr{modal}\rightarrow \mr{nodal}}(\widetilde\sigma_{mn}).
    \end{align*}
    
    \State Post-processing of $\sigma_{ij}$:
    \begin{enumerate}
        \item Guarantee positivity: $\sigma_{ij} = \sigma_{ij} - \min_{ij}\sigma_{ij}$. This is because extinction and scattering coefficients are non-negative.
        \item Normalization: $\sigma_{ij} = \frac{\sigma_{ij}}{\max_{ij}\sigma_{ij}}$.
        \item Amplification: $\sigma_{ij}=X_{A_\sigma}\sigma_{ij}$, where $X_{A_\sigma}$ is a uniform distribution on $[0,A_\sigma]$ and $A_\sigma$ is a given constant.
    \end{enumerate}
\end{algorithmic}
\end{algorithm}

% \begin{enumerate}
%     \item Draw samples of the modal representation $\widetilde\sigma_{mn}$ with low-frequency bias:
%     \begin{align*}
%         \widetilde\sigma_{mn} = \exp{\left(-c_\mr{sm}\left((\frac{m}{p_x})^2-(\frac{n}{p_y})^2\right)\right)}(X_{mn}-\frac{1}{2}),
%     \end{align*}
%     where $X_{mn}$ are i.i.d uniform distribution on $[0,1]$, and $c_\mr{sm}$ determines the decaying speed of the high-frequency signals in the distribution $\sigma_h$ so a larger $c_\mr{sm}$ is likely to generate a smoother $\sigma_h$. {This procedure is guided by the consideration that low-frequency bias exists universally in the physical world.}
    
%     \item Perform transformation from modal to nodal representations:
%     \begin{align*}
%         \sigma_{ij} = T_{\mr{modal}\rightarrow \mr{nodal}}(\widetilde\sigma_{mn}).
%     \end{align*}
    
%     \item Post-processing of $\sigma_{ij}$:
%     \begin{enumerate}
%         \item Guarantee positivity: $\sigma_{ij} = \sigma_{ij} - \min_{ij}\sigma_{ij}$. {This is because extinction and scattering coefficients are non-negative.}
%         \item Normalization: $\sigma_{ij} = \frac{\sigma_{ij}}{\max_{ij}\sigma_{ij}}$.
%         \item Amplification: $\sigma_{ij}=X_{A_\sigma}\sigma_{ij}$, where $X_{A_\sigma}$ is a uniform distribution on $[0,A_\sigma]$ {and $A_\sigma$ is a given constant.}
%     \end{enumerate}
% \end{enumerate}

We remark that in the above data generation process, the only undetermined parameters are $c_\mr{sm}$ and $A_\sigma$, which determine the overall smoothness and amplitude of the optical coefficient $\sigma_h$ which are sampled from the distribution $\mu_\sigma$, respectively.

\section{Numerical results}
\label{sec:num_exp}
In this section, we conduct numerical experiments to test the performance of our proposed element learning methods.
In the first subsection, we show how we perform the training for the neural network for element learning. Then in the following subsections, we test the performance of the trained element on a variety of numerical experiments.

\subsection{Element training}
\label{sec:num_elem_train}
As is discussed in Section \ref{sec:net_struct}, it suffices to consider the training on the reference element $\widehat K=[-1,1]^2$. We choose the (spatial) polynomial degrees $p_x=p_y=6$, namely a $Q_6$ element for the spatial discretization.
We use a fixed angular discretization of the domain $[0,2\pi]$ into $N_a=28$ uniform partitions with piece-wise constant approximation ($p_a=0$), namely a finite volume angular discretization. 

For the data generation (see Section \ref{sec:data_gen} for details), we choose the single scattering albedo $\tilde\omega=1$ and the asymmetric parameter $g_\mr{asym}=0.8$. These are standard values used for short-wave radiation of water clouds \cite{slingo1989gcm}. Following the procedures introduced in Section \ref{sec:data_gen}, we generate $N_\mr{samp}=1000$ samples of the scattering coefficient $\sigma_s$ with the amplitude parameter $A_\sigma=10$ and the smoothness parameter $c_\mr{sm}=2$. Among these samples, $800$ of them constitute the training set and the remaining $200$ samples constitute the testing set.

For the neural network structure, we construct a fully connected neural network and consider four cases for the number of layers: $N_\mr{layer}=1,2,3,4$. The size of the hidden layer is chosen to be twice of the input size, namely $N_\mr{hidden}=2N_\mr{in}$. We use the exponential linear unit (ELU) activation function for the hidden layers and linear activation for the output layer.  

The training is done on the PyTorch platform. 
We use the mean absolute error (of the neural network output, i.e., the $A_{i2o}^K$ and $A_{i2m}^K$ matrices from \eqref{eq:param2sol}) as the loss function and the Adam optimization that is innately built into the platform.
For many of these setup details for the neural network and training, we are following similar approaches as in other works such as \cite{kovachki2023neural,lu2021learning}.
We choose the batch size to be $50$. For the training, we choose the learning rate $10^{-3}$ for the first $3000$ epochs. 
We then decrease the learning rate to $10^{-4}$ and $10^{-5}$ for the following $3000$ epochs and the last $3000$ epochs, respectively.
{After the training, the NN was evaluated on Matlab by transforming into the Open Neural Network eXchange (ONNX) format.}

{The NN was trained on a single GTX30 series GPU within 20 minutes. Note that the training is required only for one time on the reference element, then this NN can be used for all the subsequent calculation without the need to retrain. }
{\color{black}
We also emphasize that it is {\sl only} the {\sl training} of the NN that is done on a GPU. Once the training is finished and the NN is transformed into an ONNX format file, only the CPU is used for the {\sl evaluation} of this NN. This is to ensure that both the element learning and the traditional numerical solvers are implemented on a CPU so that a fair comparison can be made.
}

\begin{figure}[H]
    \centering
    \includegraphics[width=0.49\textwidth]{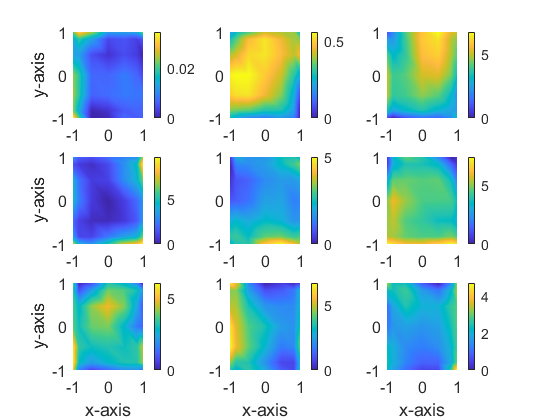}
    \includegraphics[width=0.49\textwidth]{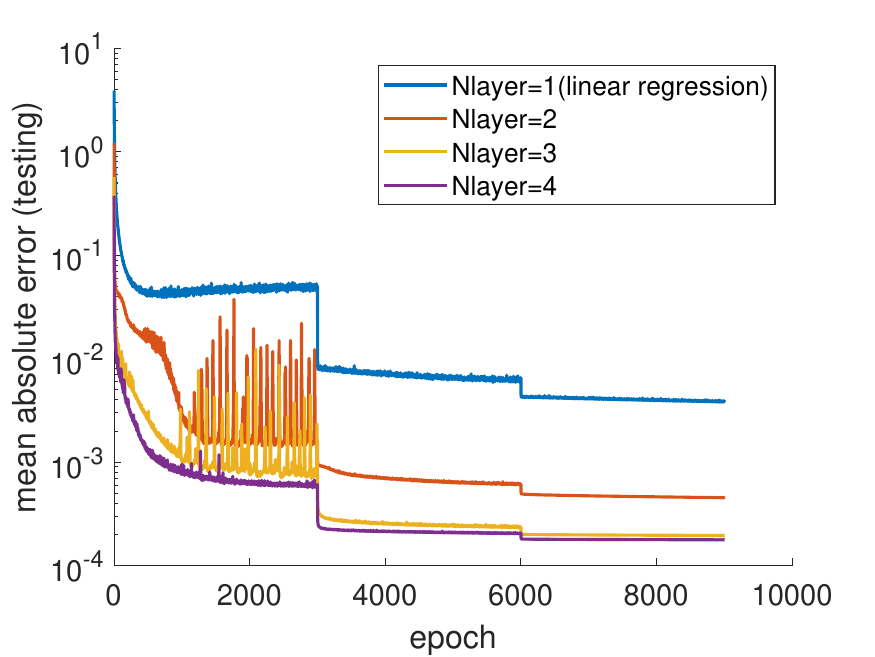}
    \caption{Training the neural network. Left: $9$ samples of the scattering coefficient $\sigma_s$ generated from the Algorithm \ref{algrm:opt_data_gen}. Right: testing error for the networks with different number of layers. The testing error is lowest (and similar) for 3 or 4 layers. The neural network is able to learn the nonlinear operator from \eqref{eq:param2sol} with an error that is nearly as small as $10^{-4}$.}
    \label{fig:samp_train_test_error}
\end{figure}

In Figure \ref{fig:samp_train_test_error}, we observe that the randomly generated scattering coefficient $\sigma_s$ presents a certain level of smoothness thanks to the low-frequency bias we introduced in the generation procedures; see Algorithm \ref{algrm:opt_data_gen}. For the approximation errors of the neural networks, we observe that the error decreases as we increase the depth of the networks. This decreasing of the error is significant when it goes from linear regression ($N_\mr{layer}=1$) to the case with one hidden layer ($N_\mr{layer}=2$). The decreasing of the error becomes marginal from $N_\mr{layer}=3$ to $N_\mr{layer}=4$.

The numerical results suggest the necessity of using neural networks with hidden layers and non-linear activation function. 
Also, the results show that a shallow network ($N_\mr{layer}=3$ or $4$) can already provide a good approximation, and in the present setup it is unlikely that any significant improvement of the approximation can be achieved by using a deeper network.

\begin{figure}[H]
    \centering
    \includegraphics[width=0.8\textwidth]{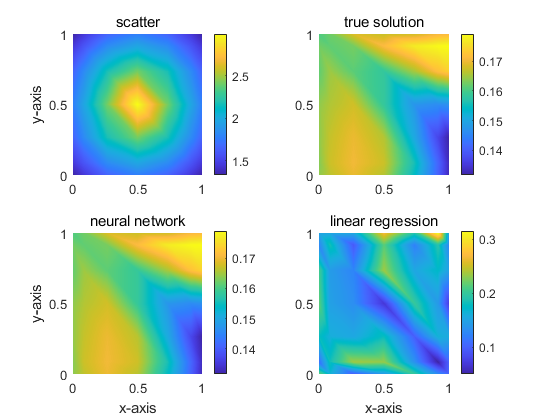}
    \caption{A comparison between the radiative transfer solution (mean intensity) obtained based on a neural network ($N_\mr{layer}=4$) versus based on linear regression ($N_\mr{layer}=1$). In this example, the inflow radiation comes from the top and left boundaries and the forcing is zero. The scattering coefficient is shown in the top-left sub-figure. The neural network provides a good approximation to the ``true solution'' (which is obtained by inverting the HDG matrix system) while the linear regression fails to provide a satisfactory approximation.}
    \label{fig:elem_lin_vs_nn}
\end{figure}

To visualize more clearly the necessity of using neural networks with hidden layers and non-linear activation, we compare the solution (mean intensity) obtained by the 4-layer network and the solution obtained by linear regression on one element. The results are collected in Figure \ref{fig:elem_lin_vs_nn}, where it is clearly shown that the neural network can provide a good approximation while the linear regression fails to produce a valid solution. 

{\color{black}
In addition to the depth of the NN, we also test how the width $N_\mr{hidden}$, the number of samples $N_\mr{samp}$, and the activation function would affect the error of the NN approximations. We also present the error with different polynomial degrees of the $Q_4$ and $Q_8$ elements. These results are collected in Figure \ref{fig:samp_train_test_error_more}.

In the top-left panel of Figure \ref{fig:samp_train_test_error_more}, we observe that by increasing the hidden layer width from $N_\mr{in}$ to $2N_\mr{in}$, there is an obvious improvement in decreasing the error, whereas the improvement from $2N_\mr{in}$ to $3N_\mr{in}$ is small. A similar pattern appears in the test with different sample numbers (top-right panel): while there is a small improvement from $N_\mr{samp}=500$ to $N_\mr{samp}=2000$, the improvement from $N_\mr{samp}=1000$ to $N_\mr{samp}=2000$ is hard to recognize. This seems to suggest that $1000$ samples is enough in the current setting of training an element.

We observe a significant improvement (more than $10$ times error reduction) by using the ELU activation function compared to the popular ReLU activation function (bottom-left panel). In the bottom-right panel, we observe that the error decreases as we increase the polynomial degree. This is likely because the network width $N_\mr{hidden}$ increases as the polynomial degree increases, since $N_\mr{hidden}=2N_\mr{in}$ and $N_\mr{in}=(p_x+1)(p_y+1)$.

\begin{figure}[H]
    \centering
    \includegraphics[width=0.49\textwidth]{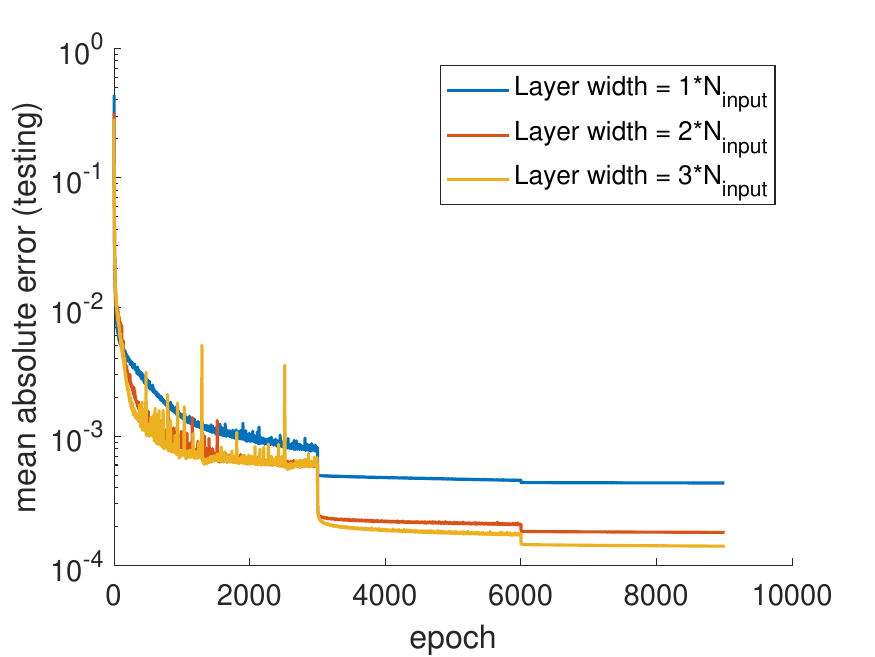}
    \includegraphics[width=0.49\textwidth]{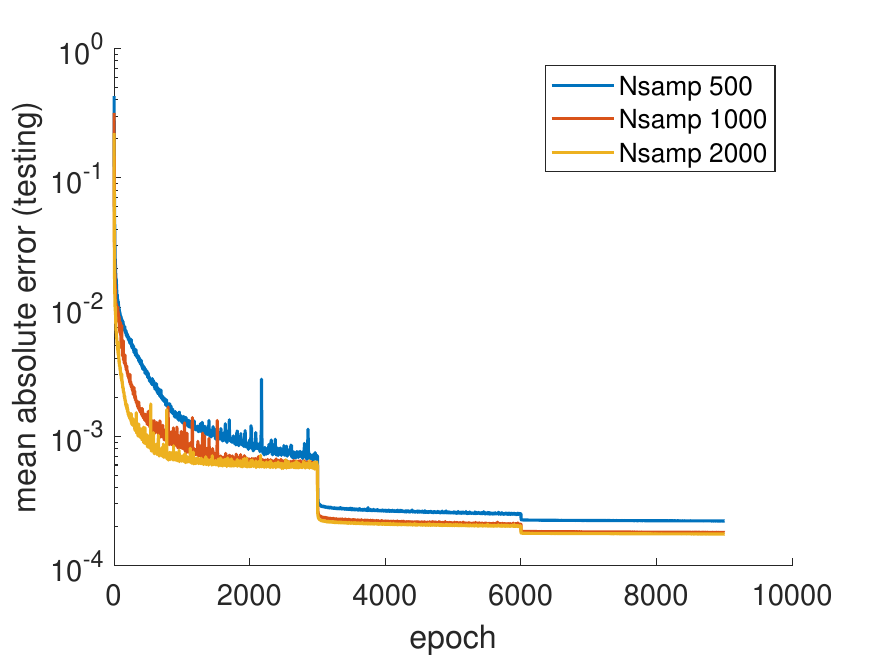}
    \includegraphics[width=0.49\textwidth]{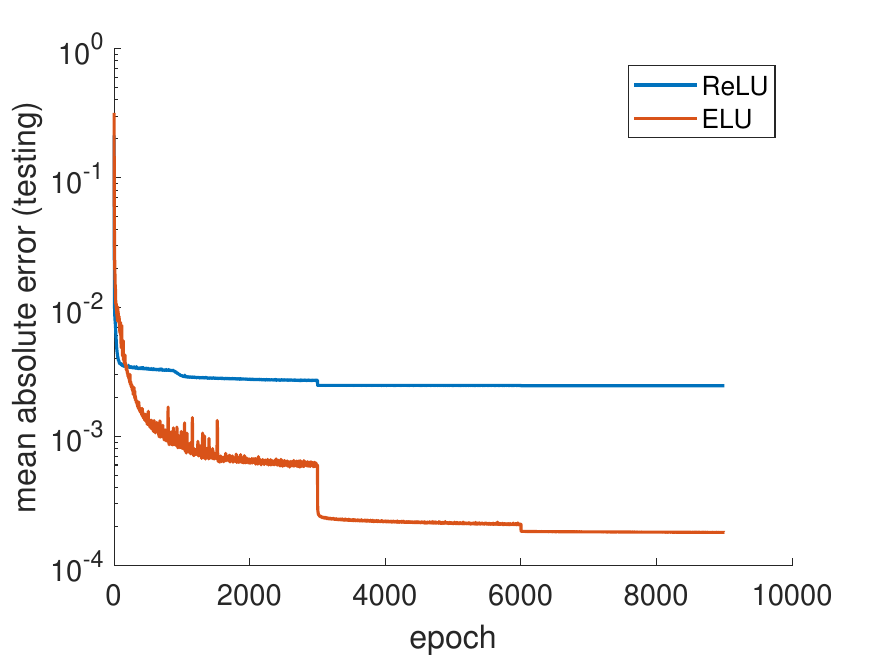}
    \includegraphics[width=0.49\textwidth]{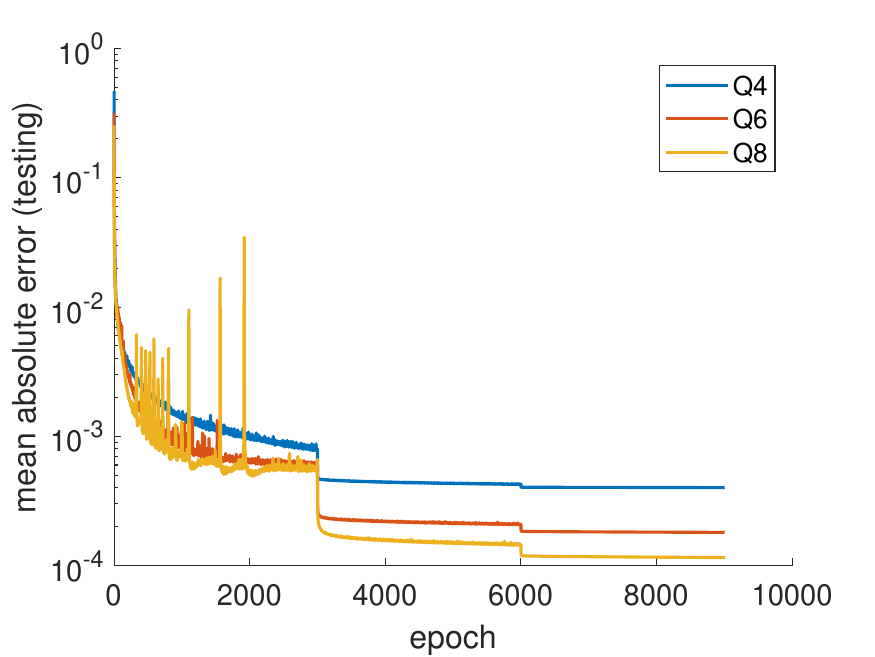}
    \caption{\color{black}
    Testing error for the networks with different hidden layer width (top-left), number of samples (top-right), activation function (bottom-left), and polynomial degree (bottom-right). These four tests are modified based on the baseline test of Figure \ref{fig:samp_train_test_error} which has the  settings of (1) $Q_6$ element, (2) number of layers $N_\mr{layer}=4$, (3) hidden layer width $N_\mr{hidden}=2N_\mr{in}$, (4) number of samples $N_\mr{samp}=1000$, and (5) the ELU activation function. Each test only changes one factor of the setting and keeps the remaining four factors unchanged.}
    \label{fig:samp_train_test_error_more}
\end{figure}

}

{\color{black}
For the rest of the numerical tests, we shall use the trained neural network with the following setting: number of layers $N_\mr{layer}=4$, hidden layer width $N_\mr{hidden}=2N_\mr{in}$, number of samples $N_\mr{samp}=1000$, and the ELU activation}.

\subsection{Testing the performance of element learning}
Here we aim to test the performance of element learning method by comparing it with other classical numerical methods. Here we consider three methods: (1) DG -- the discontinuous Galerkin method, (2) HDG -- the hybridizable discontinuous Galerkin method, and HDG-EL -- the HDG method accelerated by element learning. 

\subsubsection{Test 1 - idealized clouds}
\label{sec:test_1_ideal_cld}
\begin{figure}[H]
    \centering
    \includegraphics[width=0.49\textwidth]{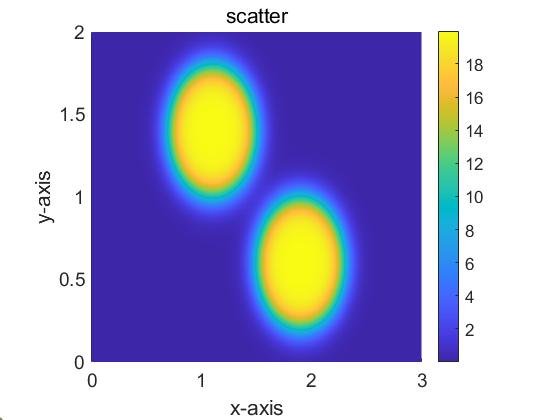}
    \includegraphics[width=0.49\textwidth]{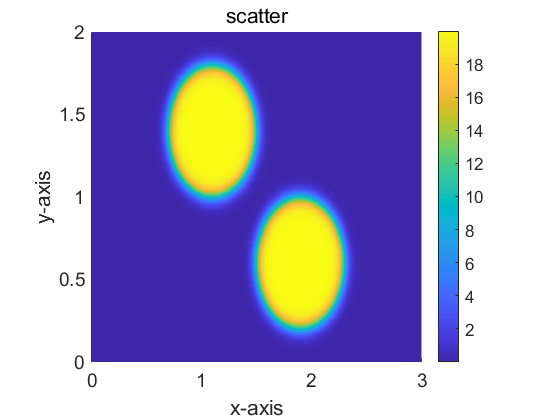}
    \caption{The distribution of the scattering coefficient $\sigma_s$ for two different test cases. Case 2 (right) has a sharper boundary transition at cloud edge compared to case 1 (left).}
    \label{fig:id_cld_sct}
\end{figure}

In this first test, we consider the case in which the radiation is scattered by idealized clouds. We consider a rectangular domain $\Omega=[0,L_x]\times[0,L_y]$ with $L_x=3$ and $L_y=2$, and two round-shaped scatterers located at the center of the domain as idealized clouds; see Figure \ref{fig:id_cld_sct}. We choose the single scattering albedo of the scatters to be $\tilde\omega=1$ and the asymmetric parameter $g_\mr{asym}=0.8$, to be consistent with the data samples we use for the element training; see Section \ref{sec:num_elem_train}.

For the boundary condition, we consider the inflow radiation $g$ to be a collimated beam coming from the top and the left sides of the domain $\Omega$, to represent solar radiation propagating towards the bottom-right direction:
\begin{align*}
    g_{\mr{left},\, \mr{top}}(\theta) =\frac{28}{2\pi} \chi_{(2\pi\frac{22}{28},2\pi\frac{23}{28})}(\theta),\qquad g_{\mr{right},\, \mr{bottom}}(\theta) =0,
\end{align*}
where $\chi$ is the characteristic function. The forcing term $f$ is chosen to be $0$.

To compare, we use the same discretization setting for the DG, HDG, and HDG-EL methods. For each element $K\in\mc T_h$, we choose $p_x=p_y=6$, $N_a=28$, and $p_a=0$ to be consistent with the setting we use for the element training; see Section \ref{sec:num_elem_train}. 
For a sequence of tests at different levels of refinement, we shall start with an initial mesh discretized by a uniform $6\times 4$ partitioning (refinement level $l=0$ with total $24$ elements). Then, for each following refinement (refinement level $l=1,2,3,4$), we discretize the domain by a uniform $3(l+2)\times 2(l+2)$ partitioning.

For the iterative solvers,
we use a stopping criterion for all the GMRES methods of $10^{-4}$ in the relative $l^2$ error. Refer to Section \ref{sec:solver} for the details of the solvers. We choose this weak stopping criterion for a fair comparison between the classical approaches (DG and HDG) and the machine learning accelerated approach (HDG-EL), since the neural network's testing error is no smaller than $10^{-4}$ (see Figure \ref{fig:samp_train_test_error}). All computation will be carried out in Matlab-R2023a with an Apple M1 CPU using a single thread.
Despite this, we remark that for HDG and HDG-EL, their local solvers could be easily parallelized with Matlab's parfor command, and a faster convergence can be expected under parallelization. Here we deliberately do not to use any parallelization for a fair comparison between DG and HDG/HDG-EL.

To estimate the error, we calculate a numerical solution by the DG method on an overrefined mesh with the refinement level $l=10$, namely, on a $36\times 24$ partition. We also increase the stopping criterion of the GMRES to $10^{-8}$ to match the high precision of this overrefined solution. We then compare the solution (mean intensity) obtained by DG, HDG, and HDG-EL methods against this overrefined solution to estimate the numerical error.

\begin{figure}[H]
    \centering
    \includegraphics[width=0.99\textwidth,trim={3cm 0 3cm 0},clip]{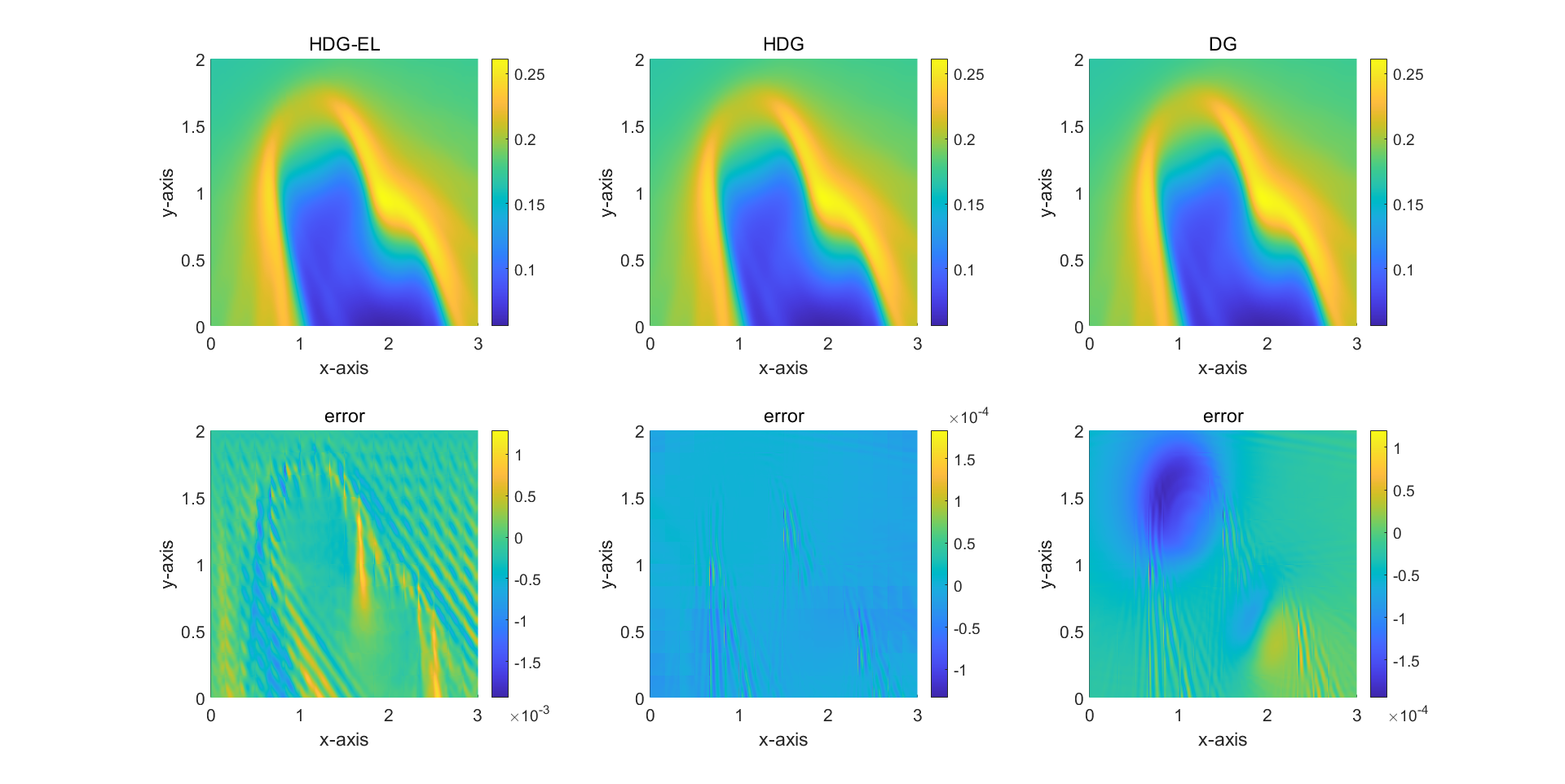}
    \caption{Numerical solution (mean intensity) and error landscape for the test case 1 (left figure of Figure \ref{fig:id_cld_sct}), with smoother transition at cloud boundary. Top row: numerical solution obtained by HDG-EL, HDG, and DG methods in the refinement level $l=4$. Bottom row: error landscape of the corresponding methods.}
    \label{fig:case1_sol_err_land}
\end{figure}

\begin{figure}[H]
    \centering
    \includegraphics[width=0.99\textwidth,trim={3cm 0 3cm 0},clip]{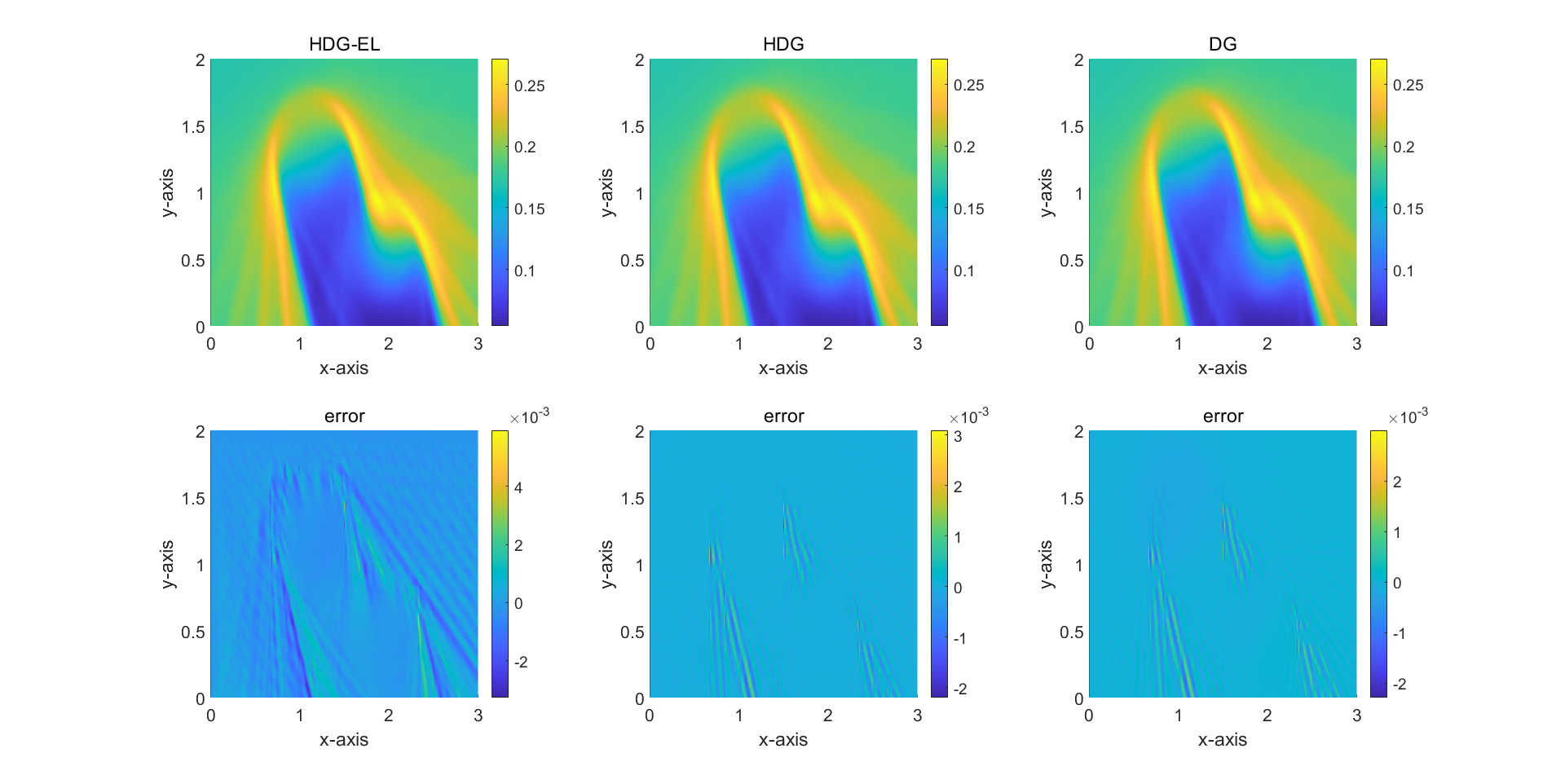}
    \caption{Numerical solution and error landscape for the test case 2 (right figure of Figure \ref{fig:id_cld_sct}), with steeper gradient at cloud boundary. Top row: numerical solution (mean intensity) obtained by HDG-EL, HDG, and DG methods in the refinement level $l=4$. Bottom row: error landscape of the corresponding methods.
    All methods achieve a similar error level of $10^{-3}$.}
    \label{fig:case2_sol_err_land}
\end{figure}

In Figure \ref{fig:case1_sol_err_land}, we plot the solution and the error obtained by HDG-EL, HDG, and DG method for the test case 1. We observe that HDG-EL provides a good approximation which reduces the error level to $10^{-3}$, but not as small as DG and HDG, which reduces the error level to $10^{-4}$.
Figure \ref{fig:case2_sol_err_land} shows the solution and the error for the test case 2, in which the scatterer has a sharper transition layer. We observe that the approximation by DG and HDG deteriorate to the error level $10^{-3}$, but the HDG-EL can still provide a similar good approximation to the error level $10^{-3}$.

{
For the following discussion, we present the solver time spent for DG, HDG, and HDG-EL methods.
For DG, it is the time spent to invert the system \eqref{eq:dg}. For HDG-local, it is the time spent to invert all the HDG local systems \eqref{eq:hdg_local_i2u},
while for HDG-EL-local, it is the time spent to infer the in2out and in2sol operators (see \eqref{eq:param2sol}) for all elements $K$ by using the network.
Finally, for HDG-global and HDG-EL-global, it is the time used to invert the global system  \eqref{eq:hdg_global_mat}.

In these comparisons, we do not account for the time required to assemble matrices.
%even though this process can constitute a large portion of the overall computational time in certain cases \cite{woopen2014comparison}
One reason is that the time it takes to assemble matrices can vary widely depending on a multitude of factors, which can complicate a comparison. For instance, the various influential factors include the software being used (e.g., Matlab based or Fortran based codes) and/or the hardware structure (stored in RAM or VRAM).
In addition, we remark that there are matrix-free approaches available for both DG and HDG methods \cite{KrKo:2019, KrKoWa2019}, which could be considered in future work.
%Finally, in many applications such as the calculation of atmospheric radiation, it involves computing radiative quantities at multiple time steps in sequence. As a result, the matrix only needs to be assembled at the initial time step. Then, subsequent steps merely require partial updates to the matrix, due to the strong correlation of atmospheric radiative quantities across time steps.
}

We show in Table \ref{tab:time_case1} (test case 1) and Table \ref{tab:time_case2} (test case 2) the solver time spent for DG, HDG, and HDG-EL methods. 
Similar results are observed in both test cases. In the last refinement level, we observe that the local solvers for the HDG method take around 36 seconds, while the global solver only takes around 0.9 second. For the HDG-EL method, the global solver takes a time similar to the HDG method. However, for the local solver, we observe a significant decreasing of the time by using element learning, which reduces the time from 36 to only 0.7 seconds, which is around $50$ times faster. The DG method only has the global system to invert, which takes around $26$ seconds, which is faster than HDG but much slower than HDG-EL.

% \begin{table}[H]
%     \centering
% \begin{tabular}{cccccc}
% \toprule
%    DOFs     &   DG    &  HDG-local & HDG-global  &  HDG-EL-local & HDG-EL-global\\
% \midrule
%   32928.00  &   4.44  &  5.60  &  0.07 &  0.90 & 0.07\\
%   74088.00  &  10.80  & 10.91  &  0.22 &  1.08 & 0.22\\     
%  131712.00  &  21.85  & 18.46  &  0.55 &  1.83 & 0.54\\   
%  205800.00  &  39.18  & 28.52  &  0.88 &  1.28 & 0.88\\     
%  296352.00  &  51.43  & 40.54  &  1.48 &  1.30 & 1.45\\     
% \bottomrule
% \end{tabular}
% \caption{Time (second) spent on DG solver, HDG local solver, HDG global solver, HDG-EL local solver, and HDG-EL global solver, for test case 1.}
% \label{tab:time_case1}
% \end{table}

\begin{table}[H]
    \centering
\begin{tabular}{cccccc}
\toprule
DOFs & DG & HDG-local & HDG-global & HDG-EL-local & HDG-EL-global\\
\midrule
   32928  &   2.24  &   4.27  &   0.05  &  0.09  &  0.04\\
   74088  &   5.36  &   9.58  &   0.13  &  0.19  &  0.13\\     
  131712  &  11.21  &  16.45  &   0.35  &  0.31  &  0.34\\     
  205800  &  19.47  &  25.33  &   0.55  &  0.46  &  0.57\\     
  296352  &  25.70  &  36.28  &   0.88  &  0.69  &  0.87\\  
\bottomrule
\end{tabular}
\caption{Time (second) spent on DG solver, HDG local solver, HDG global solver, HDG-EL local solver, and HDG-EL global solver, for test case 1.}
\label{tab:time_case1}
\end{table}

% \begin{table}[H]
%     \centering
% \begin{tabular}{cccccc}
% \toprule
%    DOFs      &  DG    &  HDG-local  &  HDG-global  &  HDG-EL-local  &  HDG-EL-global\\
% \midrule
%    32928.00  &   4.63 &  4.61  &  0.08 & 0.35  & 0.08\\     
%    74088.00  &  11.30 & 10.08  &  0.22 & 0.83  & 0.21\\     
%   131712.00  &  23.37 & 17.80  &  0.54 & 0.75  & 0.58\\     
%   205800.00  &  41.97 & 27.07  &  0.88 & 0.89  & 0.89\\     
%   296352.00  &  54.73 & 39.80  &  1.49 & 1.80  & 1.49\\   
% \bottomrule
% \end{tabular}
% \caption{Time (second) spent on DG solver, HDG local solver, HDG global solver, HDG-EL local solver, and HDG-EL global solver, for test case 2.}
% \label{tab:time_case2}
% \end{table}

\begin{table}[H]
    \centering
\begin{tabular}{cccccc}
\toprule
  DOFs  & DG & HDG-local  &  HDG-global  &  HDG-EL-local  &  HDG-EL-global\\
\midrule
  32928  &   2.32  &   4.29  &  0.04  &  0.09  &  0.04\\     
  74088  &   5.57  &   9.52  &  0.14  &  0.19  &  0.13\\     
 131712  &  11.63  &  16.15  &  0.35  &  0.31  &  0.34\\     
 205800  &  20.71  &  25.33  &  0.56  &  0.48  &  0.54\\     
 296352  &  26.77  &  36.60  &  0.90  &  0.70  &  0.89\\     
\bottomrule
\end{tabular}
\caption{Time (second) spent on DG solver, HDG local solver, HDG global solver, HDG-EL local solver, and HDG-EL global solver, for test case 2.}
\label{tab:time_case2}
\end{table}

\begin{figure}[H]
    \centering
    \includegraphics[width=0.32\textwidth]{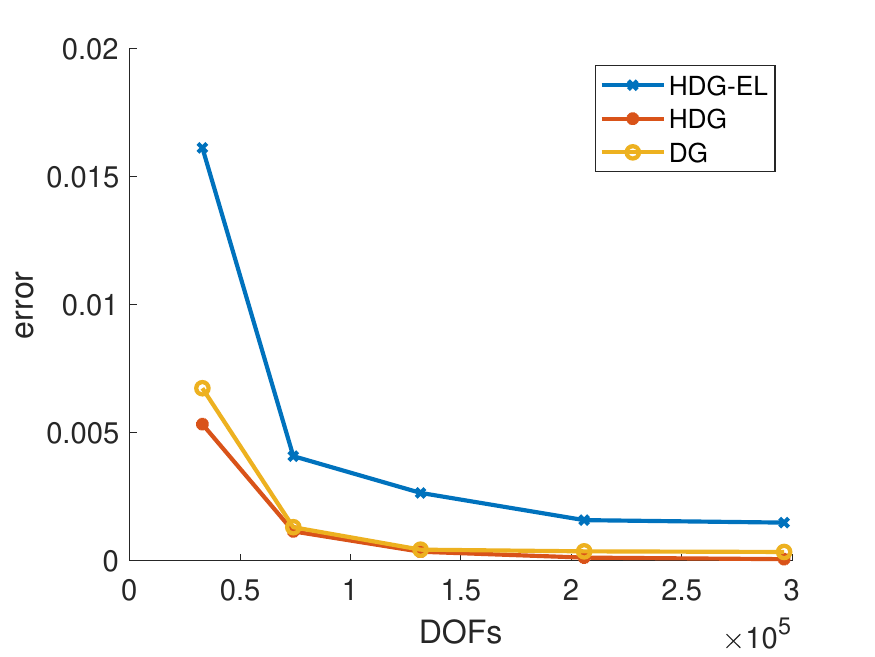}
    \includegraphics[width=0.32\textwidth]{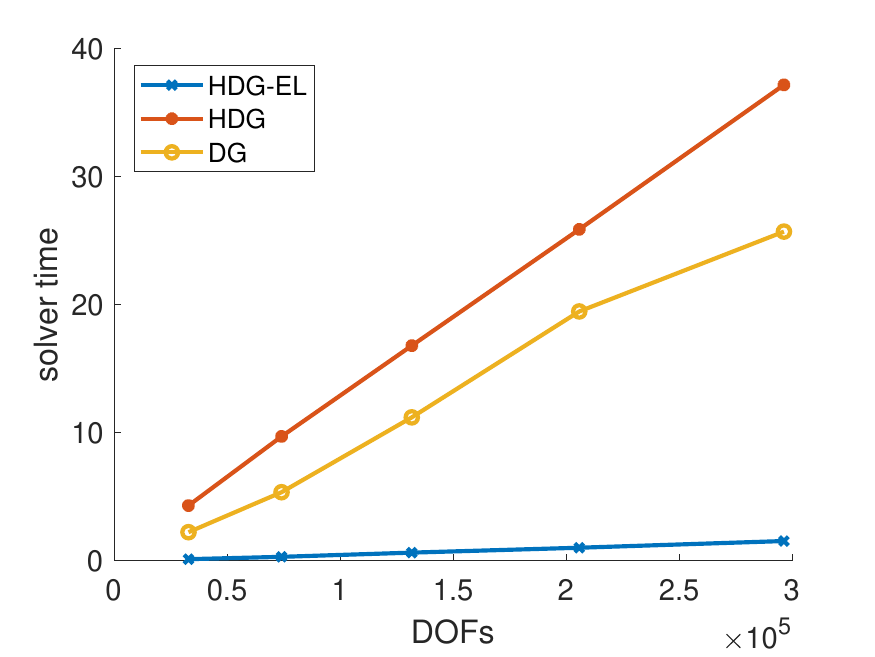}
    \includegraphics[width=0.32\textwidth]{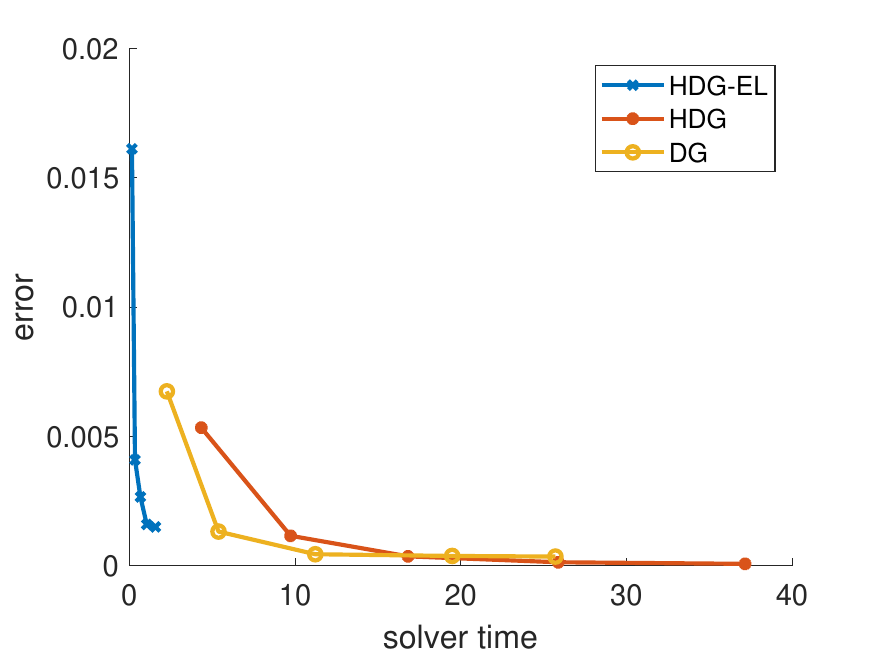}
    \caption{Comparison of the different methods in terms of cost and accuracy, for test case 1. Left: DOFs vs relative $L^2$ error. Middle: DOFs vs solver time. Right: solver time vs relative $L^2$ error. The solver time for HDG and HDG-EL is the summation of the local solver and the global solver time.  
    For an error level of approximately $2\times 10^{-3}$, in the right panel, the HDG-EL method is around $5$ to $10$ times faster than the DG and HDG methods.}
    \label{fig:dof_err_time_case1}
\end{figure}

\begin{figure}[H]
    \centering
    \includegraphics[width=0.32\textwidth]{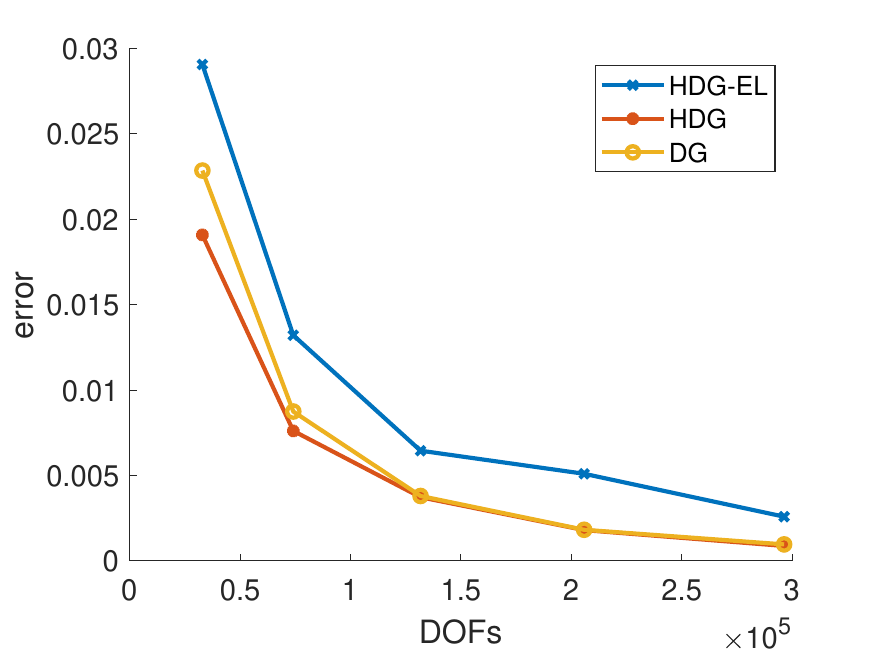}
    \includegraphics[width=0.32\textwidth]{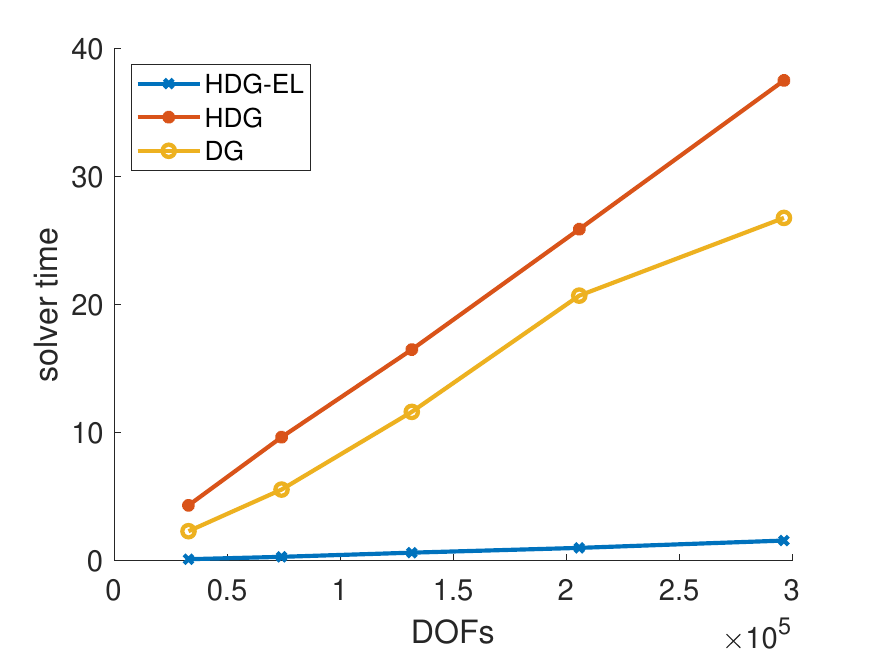}
    \includegraphics[width=0.32\textwidth]{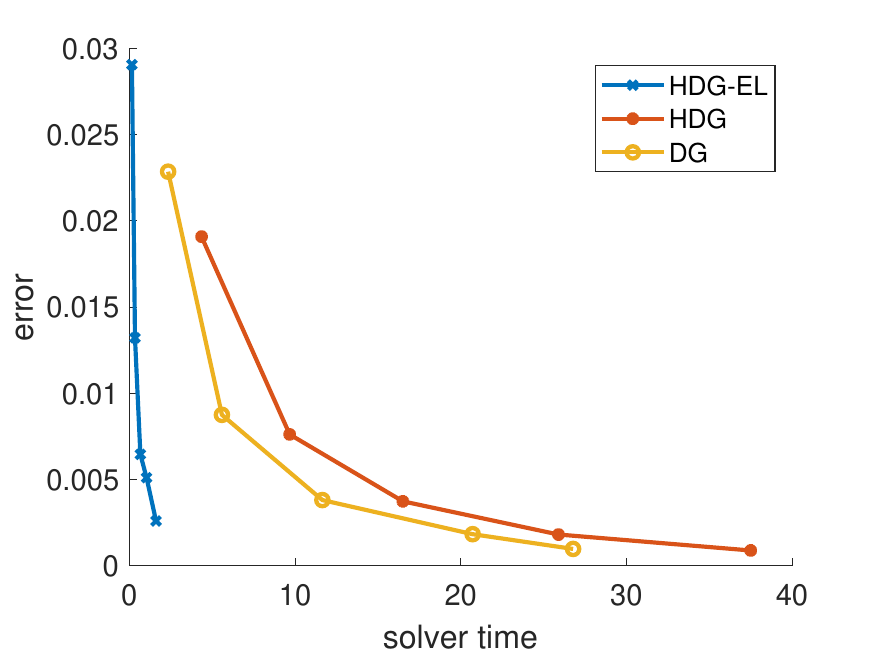}
    \caption{Comparison of the different methods in terms of cost and accuracy, for test case 2. Left: DOFs vs relative $L^2$ error. Middle: DOFs vs solver time. Right: solver time vs relative $L^2$ error.  The solver time for HDG and HDG-EL is the summation of the local solver and the global solver time.
    For an error level of approximately $2\times 10^{-3}$, in the right panel, the HDG-EL method is around $10$ times faster than the HDG and DG methods.}
    \label{fig:dof_err_time_case2}
\end{figure}

To better see the relations between the DOFs, the error, and the solver time for the three methods, we plot Figure \ref{fig:dof_err_time_case1} and Figure \ref{fig:dof_err_time_case2} for case 1 and case 2, respectively. 
By the left sub-figures of Figure \ref{fig:dof_err_time_case1} and Figure \ref{fig:dof_err_time_case2}, we observe that all the three methods can reduce the relative $L^2$ error to a satisfactory level of $10^{-3}$. In addition, we observe that HDG and DG are more efficient in reducing the error (with respect to the DOFs) compared to the HDG-EL method. This is expected since HDG-EL is an approximation to the HDG method by using neural networks. 
We also observe that, for the test case 1 (which has a smoother transition layer for the scatterer; see Figure \ref{fig:id_cld_sct}), the errors are smaller for all the three methods, compared to the test case 2, which has a steeper transition layer for the scatterer.

In the middle figures of Figure \ref{fig:dof_err_time_case1} and Figure \ref{fig:dof_err_time_case2}, we plot how the solver time changes according to the increase of the DOFs. For both test case 1 and test case 2, we observe that HDG-EL is significantly faster than HDG and DG methods, while the DG method is faster than the HDG method.

In the right figures of Figure \ref{fig:dof_err_time_case1} and Figure \ref{fig:dof_err_time_case2}, we plot how the error changes according to the solver time. We observe that the HDG-EL methods are much faster than HDG and DG methods with a fixed error level. In the error level around $2\times 10^{-3}$, we observe the HDG-EL method is about $5$ to $10$ times faster than the DG and HDG methods (around 5 to 10 times faster for test case 1, and around 10 times faster for test case 2). 
It is interesting to point out that despite the overall faster convergence of the HDG-EL method compared to HDG/DG in both the test cases, the speed-up of the HDG-EL method is more obvious in a more difficult test case (test case 2 with sharper transition at scatterer boundary). This seems to suggest the NN by element learning generalizes well (to relatively steep gradients in the PDE parameters and solutions).

{\color{black}
We next show the numerical tests for $Q_4$ and $Q_8$ elements, to be compared with the results of the $Q_6$ element in Table \ref{tab:time_case2} and Figure \ref{fig:dof_err_time_case2}. All the test settings will be kept the same except that the polynomial degree will be changed.

Table \ref{tab:time_case2_Q4} and Figure \ref{fig:dof_err_time_case2_Q4} show the test results for  $Q_4$ element, while 
Table \ref{tab:time_case2_Q8} and Figure \ref{fig:dof_err_time_case2_Q8} show the results for  $Q_8$ element. 
By comparing the left panel of the Figure \ref{fig:dof_err_time_case2_Q4} and Figure \ref{fig:dof_err_time_case2_Q8},
we observe that HDG-EL is similar to HDG and DG in accuracy with fixed DOFs for $Q_4$ element. However, HDG-EL is much less accurate than DG and HDG with fixed DOFs for $Q_8$ element.
An explanation for this observation is that the traditional high-order finite element solvers become highly accurate with only a few element. However, the element learning approach introduces machine learning error (or neural network approximation error) which does not decrease as fast as the high-order finite element methods. 

When we compare the solver time while keeping the DOFs fixed (middle panel of Figure \ref{fig:dof_err_time_case2_Q4} and Figure \ref{fig:dof_err_time_case2_Q8}), we observe that the speed-up from element learning is more significant in the case of the higher order $Q_8$ methods, compared to the lower order $Q_4$ methods. This is most obvious by comparing the local solver time of HDG and HDG-EL; see Table \ref{tab:time_case2_Q4} and Table \ref{tab:time_case2_Q8}. While the speed-up varies from $7$ to $9$ for $Q_4$ element, the speed-up for $Q_8$ element varies from $25$ to $65$.

Finally, for the solver time versus the error (right panel of Figure \ref{fig:dof_err_time_case2_Q4} and Figure \ref{fig:dof_err_time_case2_Q8}), 
we again observe an approximately $5$ to $10$ times speed-up similar to the case of $Q_6$.
Also note that the error in element learning does not reach below
$O(10^{-3})$, due to the machine learning component of element learning.
For even higher accuracy with errors below $O(10^{-3})$,
the corresponding machine learning part needs to provide a more precise approximation.
The present results here are able to achieve low errors of
$O(10^{-1})$ to $O(10^{-3})$ which are sufficiently accurate for many applications.
% Notice that while element learning for the $Q_8$ methods can reduce the solver time significantly for fixed DOFs, it is not sufficiently accurate in approximating the FEM solver. 
% This seems to indicate that to fully realize the potential of high-order element learning, the corresponding machine learning part needs to provide a more precise approximation.

\begin{table}[H]
    \centering
\begin{tabular}{cccccc}
\toprule
  DOFs  & DG & HDG-local  &  HDG-global  &  HDG-EL-local  &  HDG-EL-global\\
\midrule
 16800.00 & 0.83 & 0.81  &  0.04 &  0.07  &  0.03\\
 37800.00 & 2.10 & 1.76  &  0.08 &  0.19  &  0.08\\     
 67200.00 & 3.86 & 3.04  &  0.20 &  0.24  &  0.19\\     
105000.00 & 6.06 & 5.06  &  0.31 &  0.39  &  0.31\\     
151200.00 & 9.03 & 7.16  &  0.46 &  0.55  &  0.48\\    
\bottomrule
\end{tabular}
\caption{\color{black}Time (second) spent on DG solver, HDG local solver, HDG global solver, HDG-EL local solver, and HDG-EL global solver, for test case 2, for $Q_4$ element.}
\label{tab:time_case2_Q4}
\end{table}

\begin{figure}[H]
    \centering
    \includegraphics[width=0.32\textwidth]{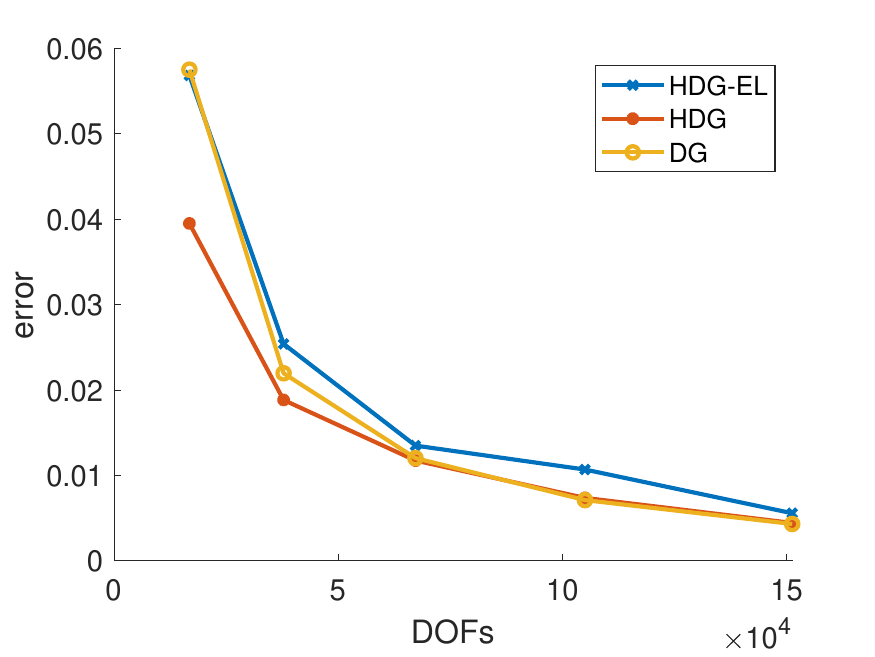}
    \includegraphics[width=0.32\textwidth]{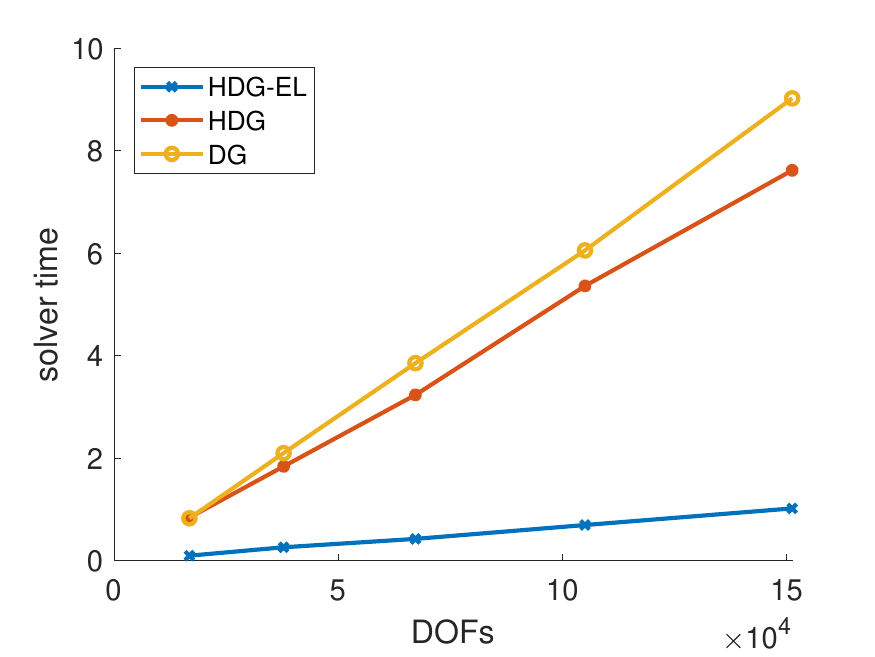}
    \includegraphics[width=0.32\textwidth]{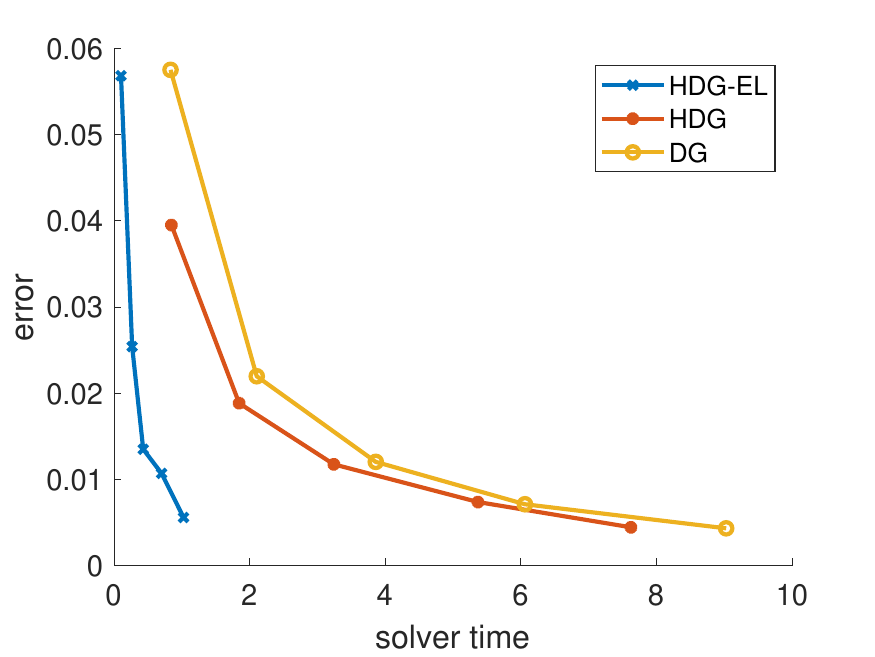}
    \caption{\color{black}Same test setting of Figure \ref{fig:dof_err_time_case2} but with $Q_4$ element.}
    \label{fig:dof_err_time_case2_Q4}
\end{figure}

\begin{table}[H]
    \centering
\begin{tabular}{cccccc}
\toprule
  DOFs  & DG & HDG-local  &  HDG-global  &  HDG-EL-local  &  HDG-EL-global\\
\midrule
54432.00  &  5.13 &  18.95 &  0.07  &  0.13 &  0.06\\     
122472.00 & 15.08 &  40.85 &  0.20  &  0.26 &  0.21\\   
217728.00 & 27.61 &  71.93 &  0.53  &  0.45 &  0.53\\     
340200.00 & 43.06 & 110.43 &  0.90  &  0.84 &  0.89\\     
489888.00 & 64.45 & 159.01 &  1.41  &  0.98 &  1.40\\     
\bottomrule
\end{tabular}
\caption{\color{black}Time (second) spent on DG solver, HDG local solver, HDG global solver, HDG-EL local solver, and HDG-EL global solver, for test case 2, for $Q_8$ element.}
\label{tab:time_case2_Q8}
\end{table}

\begin{figure}[H]
    \centering
    \includegraphics[width=0.32\textwidth]{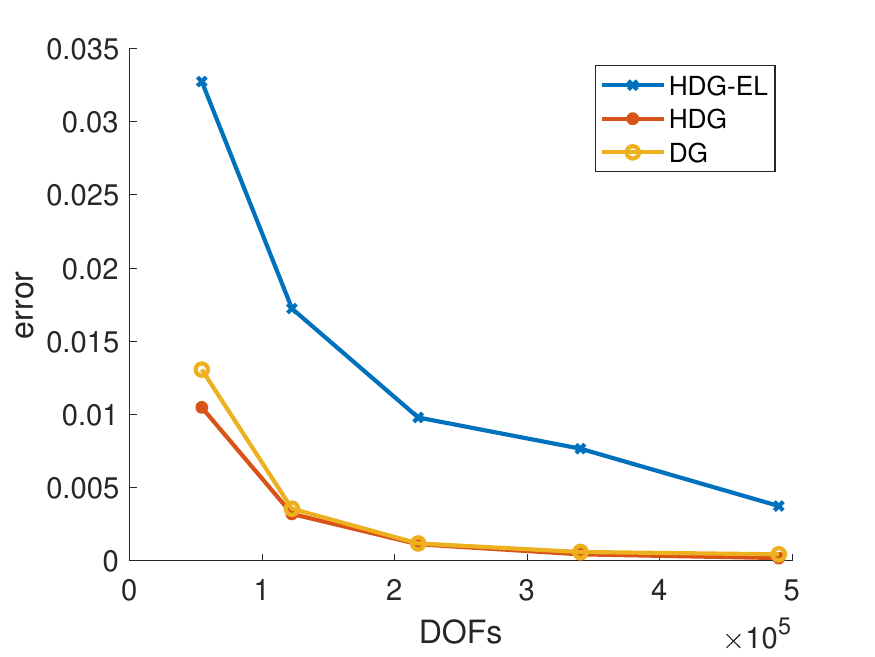}
    \includegraphics[width=0.32\textwidth]{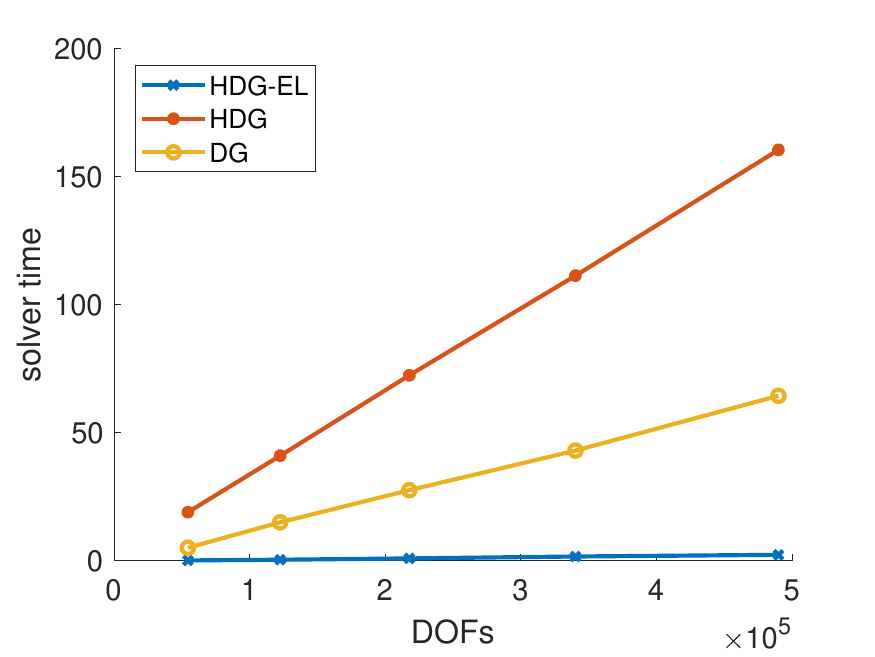}
    \includegraphics[width=0.32\textwidth]{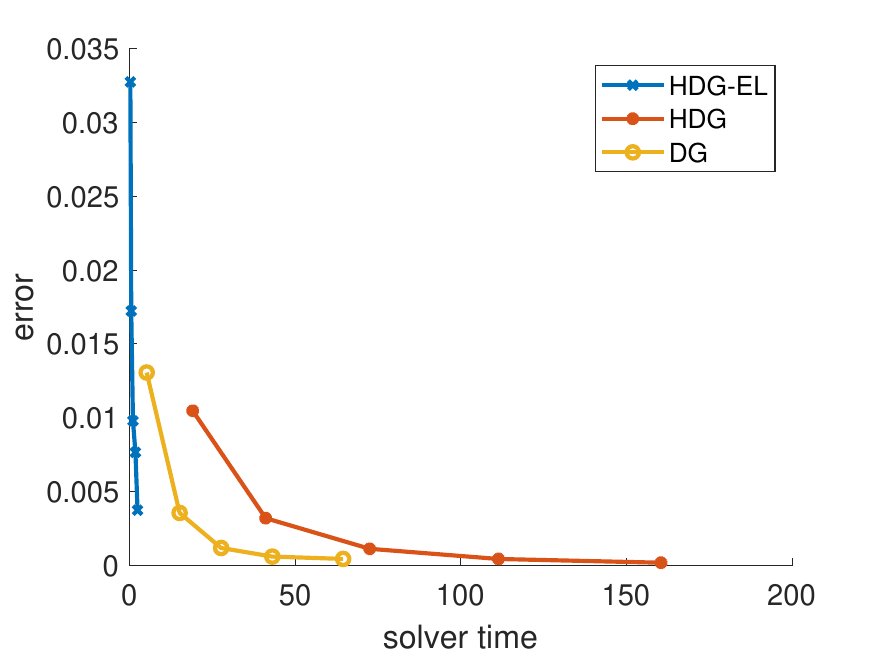}
    \caption{\color{black}Same test setting of Figure \ref{fig:dof_err_time_case2} but with $Q_8$ element.}
    \label{fig:dof_err_time_case2_Q8}
\end{figure}

}

\subsubsection{Test 2 - I3RC cloud fields}
\label{sec:num_exp_i3rc}
\begin{figure}[H]
    \centering
    \includegraphics[width=0.99\textwidth,trim={4cm 0 4cm 0},clip]{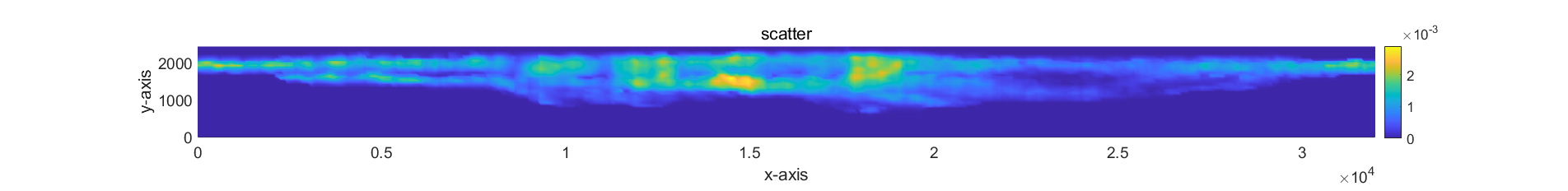}
    \caption{Cloud field from I3RC case 2.}
    \label{fig:i3rc_cloud_case2}
\end{figure}

Here we consider a more realistic test case, for which the cloud field is retrieved from the International Comparison of 3-Dimensional Radiative Transfer Codes (I3RC) case 2 \cite{cahalan2005i3rc}; see Figure \ref{fig:i3rc_cloud_case2} for a visualization of the cloud field. We shall use the same experiment setting of the tests in Section \ref{sec:test_1_ideal_cld}, unless otherwise specified.

For the boundary condition, we take the inflow radiation $g$ as follows:
\begin{align*}
    g_{\mr{left},\, \mr{top}}(\theta) =\frac{28}{2\pi} \chi_{(2\pi\frac{24}{28},2\pi\frac{25}{28})}(\theta),\qquad g_{\mr{right},\, \mr{bottom}}(\theta) =0,
\end{align*}
which mimics the solar radiation coming from the top and the left boundary and going towards the bottom-right direction. 

To compare, we use the same discretization setting for DG, HDG, and HDG-EL methods. We shall start with an initial mesh discretized by a uniform $26\times 2$ partitioning (refinement level $l=0$ with total $52$ elements). Then, for each refinement (refinement level $l=1,2,3,4$), we discretize the domain by a uniform $13(l+2)\times 1(l+2)$ partitioning.
To estimate the error, we calculate a numerical solution by the DG method on an overrefined mesh with the refinement level $l=10$, namely, on a $156\times 12$ mesh. We then compare the solution (mean intensity) obtained by HDG-EL, HDG, and DG methods against this over-refined solution to estimate the numerical error.

\begin{figure}[H]
    \centering
    \includegraphics[width=0.99\textwidth,trim={0 2cm 0 1cm},clip]{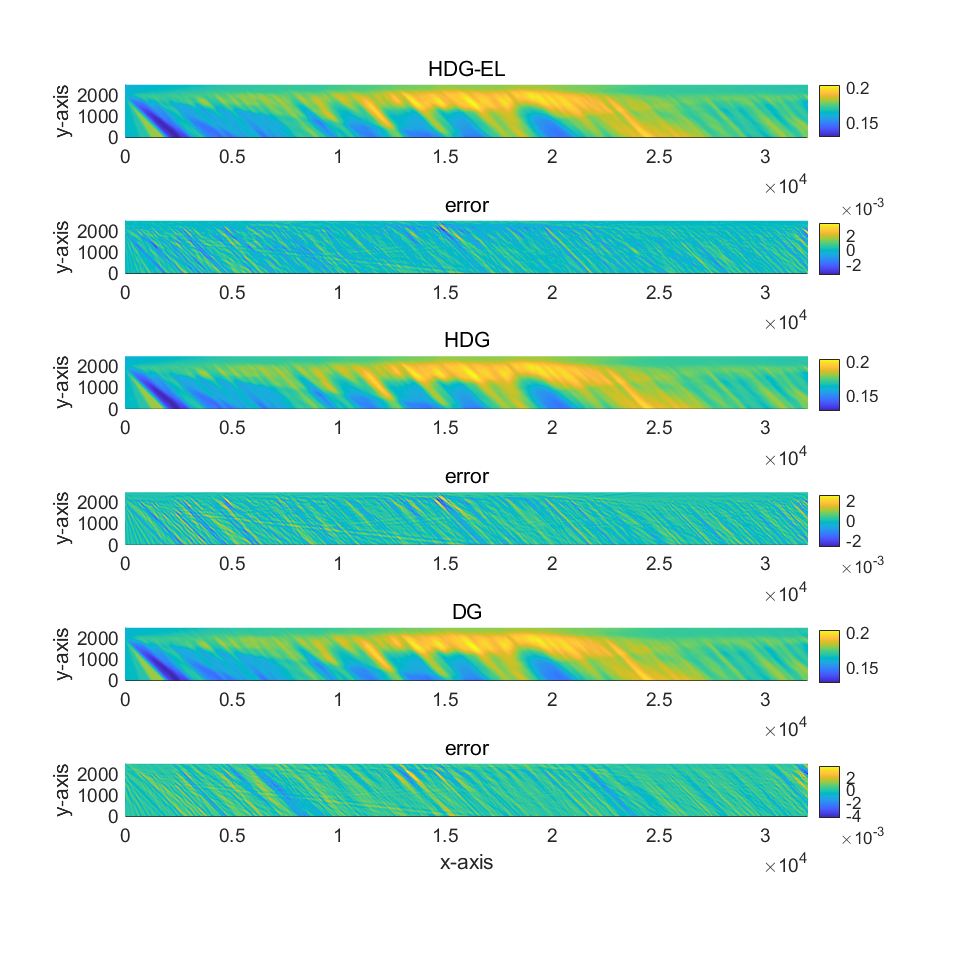}
    \caption{Numerical solution and error landscape for I3RC test case (see Figure \ref{fig:i3rc_cloud_case2}). Odd row: numerical solution (mean intensity) obtained by HDG-EL, HDG, and DG methods in the refinement level $l=4$. Even row: error landscape of the corresponding methods.
    All methods achieve a similar error level of approximately $2\times 10^{-3}$.}
    \label{fig:i3rc_sol_err_land}
\end{figure}

In Figure \ref{fig:i3rc_sol_err_land}, we plot the solution and the error obtained by the HDG-EL, HDG, and DG methods. We observe that HDG-EL provides a  good approximation to the overrefined solution. The HDG-EL method reduces the error level to $2\times 10^{-3}$, which is very close to HDG and is even slightly better than DG.

We show the solver time for the three methods in Table \ref{tab:time_i3rc}.
We again observe that HDG-EL uses much less time (around 20 times faster) for the local solver compared to the HDG method. We observe slightly slower global solver with HDG-EL compared to HDG method. Again, we observe that the DG method is faster than HDG but much slower than HDG-EL. 

% \begin{table}[H]
%     \centering
% \begin{tabular}{cccccc}
% \toprule
%    DOFs    &    DG &  HDG-local  &  HDG-global  &  HDG-EL-local  &  HDG-EL-global\\
% \midrule
%   71344.00 &  3.32 &  3.87  & 0.17 &  0.46  &  0.20\\
%  160524.00 &  9.92 & 10.16  & 0.52 &  0.63  &  0.64\\     
%  285376.00 & 21.54 & 17.82  & 1.13 &  1.08  &  1.53\\     
%  445900.00 & 36.07 & 33.98  & 1.94 &  1.62  &  2.78\\     
%  642096.00 & 55.29 & 46.26  & 3.36 &  2.46  &  5.03\\
% \bottomrule
% \end{tabular}
% \caption{Time (second) spent on DG solver, HDG local solver, HDG global solver, HDG-EL local solver, HDG-EL global solver, for the I3RC test case.}
% \label{tab:time_i3rc}
% \end{table}

\begin{table}[H]
    \centering
\begin{tabular}{cccccc}
\toprule
 DOFs & DG & HDG-local & HDG-global & HDG-EL-local & HDG-EL-global\\
\midrule
   71344  &   1.62  &   2.74  &  0.12  &  0.19 & 0.12\\
  160524  &   5.18  &   6.64  &  0.33  &  0.37 & 0.41\\     
  285376  &  11.02  &  13.27  &  0.62  &  0.66 & 0.89\\     
  445900  &  15.86  &  21.31  &  1.11  &  0.99 & 1.52\\     
  642096  &  24.73  &  33.47  &  2.01  &  1.43 & 2.92\\  
\bottomrule
\end{tabular}
\caption{Time (second) spent on DG solver, HDG local solver, HDG global solver, HDG-EL local solver, HDG-EL global solver, for the I3RC test case.}
\label{tab:time_i3rc}
\end{table}

\begin{figure}[H]
    \centering
    \includegraphics[width=0.32\textwidth]{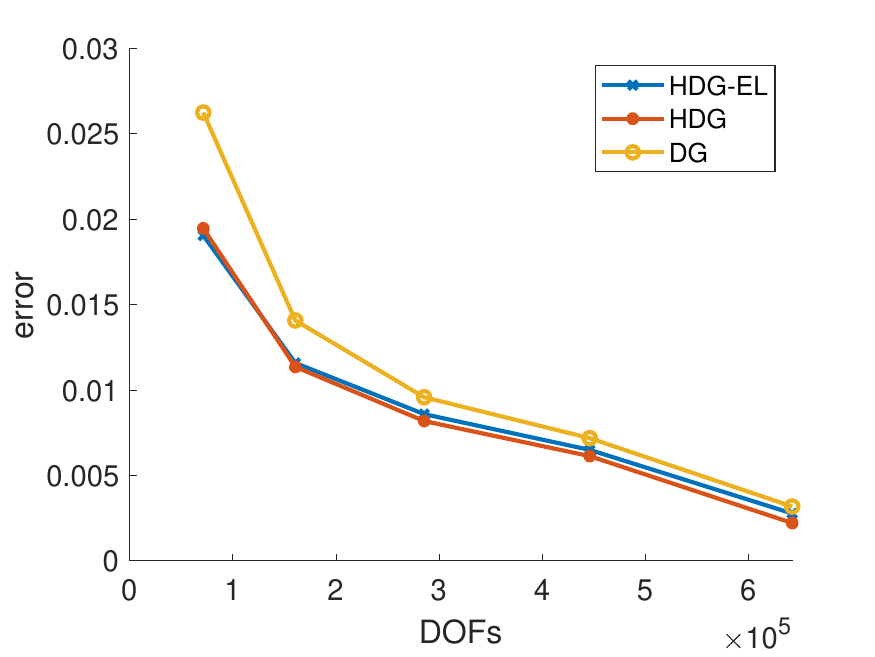}
    \includegraphics[width=0.32\textwidth]{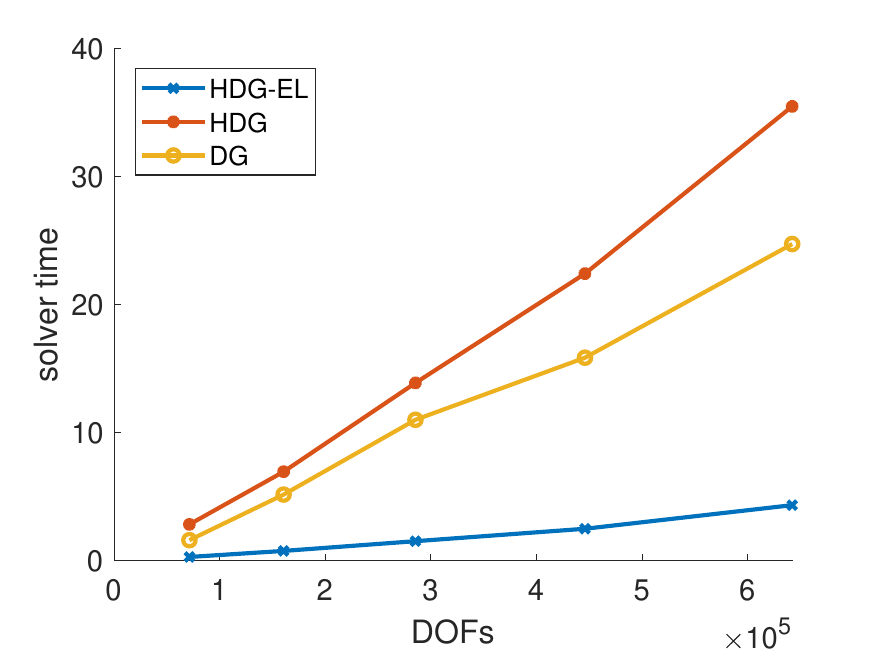}
    \includegraphics[width=0.32\textwidth]{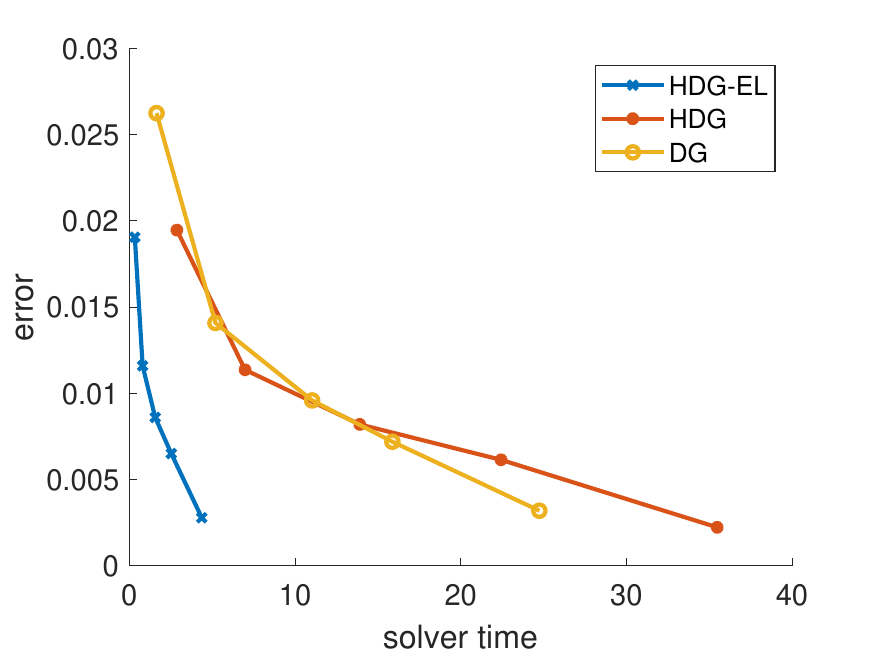}
    \caption{Comparison of the different methods in terms of cost and accuracy, for the I3RC test case. Left: DOFs vs relative $L^2$ error. Middle: DOFs vs solver time. Right: solver time vs relative $L^2$ error.  The solver time for HDG and HDG-EL is the summation of the local solver and the global solver time. 
    For an error level of approximately $2\times 10^{-3}$, in the right panel, the HDG-EL method is $5$ to $10$ times faster than the DG and HDG methods.}
    \label{fig:dof_err_time_i3rc}
\end{figure}

We again show the relations between the DOFs, the error, and the solver time in Figure \ref{fig:dof_err_time_i3rc}.
In the left sub-figure Figure \ref{fig:dof_err_time_i3rc}, we observe that all three methods behave similarly in reducing the error with respect to the DOFs. Unlike the test cases in Section \ref{sec:test_1_ideal_cld}, here we observe a similar performance of HDG-EL compared to DG and HDG. 

In the middle sub-figures of Figure \ref{fig:dof_err_time_i3rc}, we again observe a significant reduction of the solver time of HDG-EL compared to HDG and DG methods.
Because of this, in the right sub-figures of Figure \ref{fig:dof_err_time_i3rc}, we observe that HDG-EL is much faster in reducing the error to the level of $2\times 10^{-3}$, compared to HDG  and DG methods. We again observe a $5$ to $10$ times speed-up by using element learning.

\section{Conclusion}
\label{sec:conclusion}
In this paper, we propose a novel approach to accelerating finite element-type methods by machine learning. A main goal is for the numerous benefits of traditional finite element-type methods to be retained. Many of the advantageous features of the element learning approach were listed near the end of section~\ref{sec:intro}. Numerical experiments have demonstrated promising results that an approximately $5$-$10$ times speed-up can be achieved by using the proposed element learning method. 

As was pointed out above, the element learning method is closely related to hybridizable discontinuous Galerkin methods in the sense that the local solvers are replaced by machine learning approaches. As a result, all hybridizable discontinuous Galerkin methods, hybridized mixed method, and hybridized continuous Galerkin methods, can be potentially accelerated by element learning. This generality suggests the potential of the proposed method to be applied to more PDEs that have been studied with HDG discretizations. Examples include but are not limited to Maxwell's equations, linear elasticity, convection-diffusion, and Stokes/Navier-Stokes equations.

Despite the great potential of the proposed method, there still exist many topics that remains to be studied. Here we list a few of them:
\begin{enumerate}
\item Even greater speed-up can be expected if a higher order polynomial degree $p$ is used. At the same time, for larger $p$, other factors may arise, as discussed in the next item.
\item As the polynomial degree $p$ increases, a larger ratio of the DOFs can be saved by the techniques of hybridization and/or static condensation. So, a greater acceleration can be expected by the element learning approach. However, as we increase the polynomial degree, it becomes more difficult to train the neural network. This can leads to less reliable predictions by the network. Therefore, how to choose the polynomial degree for a good balance between the efficiency of the method (for which we aim to increase $p$), and the reliability of the method (for which we aim to decrease $p$), remains to be studied. 

\item The current approach achieves a computational speed-up even with using a {simple} fully connected neural network without any a priori structures. It remains to be studied that whether better approximation of the networks {or reduced computation cost} can be achieved by introducing additional structures, such as {dropout} \cite{srivastava2014dropout}, convolutional layers (CNN) \cite{lecun1998gradient}, attention mechanism (transformer) \cite{vaswani2017attention}, or shortcut connections (ResNet) \cite{he2016deep}. 

\item The testing error of the neural network for element learning did not reach below $10^{-4}$. This limitation hampers the method to be used for obtaining solutions with high accuracy requirement. It remains to be explored if additional techniques can reduce the error to a lower value.

\item In some applications, a desirable property of  data-driven approaches is the ability to learn a data-driven solver without knowing the underlying PDE or dynamical equations.  Here, for element learning, note that the underlying PDE does not need to be known. The global solver is defined by the consistency condition in \eqref{eq:hdg_global}, which makes no explicit reference to the PDE itself, and instead refers only to a local solver $u_h(\widehat u_h,f)$ in abstract form. In the HDG method, the local solver comes from the PDE itself. In element learning, on the other hand, the local solver is learned from data, which could come from synthetic data from a numerical PDE solution if the PDE is known, or could come from observational measurements or laboratory experiments. 
Knowledge of boundary conditions is still needed for an appropriate definition of the local in2out operators in element learning, in order to ensure a well-posed problem, and to allow the coupling of local solvers in creating the global solver.
%Also note that the training in element learning is required only for a single element, which is a smaller and more feasible training problem than training for the solution over the entire domain. 
It would be interesting in the future to investigate the ideas of element learning in scenarios where the underlying PDE or dynamical equations are unknown due to the complexity of the system, or are computationally intractable due to the high dimensionality of the true system. For example, such a situation arises in physics parameterizations or subgrid-scale parameterizations in atmospheric, oceanic, and climate dynamics and other complex systems \cite{rasp2018deep,wu2020enforcing,clark2022correcting,guan2022stable,boral2023neural,perezhogin2023generative}.

{\color{black}
\item It is possible to discretize the in2out and the in2sol operators on each element $K$ by numerical approaches other than the HDG local solvers. 
Examples include low-order finite element methods (macro-element) and finite difference methods. In this case, each element $K$ is more like a small domain instead of an element in the traditional sense.

In the case of macro-element, each macro-element $K$ is discretized by several low-order elements. This setting can be related to the setting of multi-scale finite element methods; see, for instance, the work \cite{EfLaSh:2015} where a multiscale HDG method was developed for second order elliptic equations. Considering the close connection between element learning and HDG methods, it would be interesting to explore the potential of element learning for solving multiscale problems.
}

\item The main idea of element learning is to approximate the map from element geometry and coefficients to the element-wise solution operator. A topic that remains to be studied is the theoretical aspects of the approximation abilities of neural networks {\color{black}\cite{ShYaZh:2022,EMaWu:2022,HeLiXuZh:2020,chen2023deep,lanthaler2023curse}} for this map. This includes three parts: (1) best approximation capability of the network, (2) the estimate of the optimization error, and (3) the estimate of the generalization error.

\end{enumerate}

% \section*{CRediT authorship contribution statement}
% {\bf Shukai Du:} Conceptualization, Methodology, Writing– Original draft preparation, Writing--Review \& Editing, Software.\\
% {\bf Samuel N. Stechmann:} Conceptualization, Methodology,  Writing--Review \& Editing, Supervision.

% \section*{Declaration of competing interest}
% The authors declare that they have no known competing financial interests or personal relationships that could have appeared to influence the work reported in this paper.

\section*{Data availability}
The cloud data used for the numerical experiment Section \ref{sec:num_exp_i3rc} can be publicly accessed through the link
\url{https://earth.gsfc.nasa.gov/climate/model/i3rc/testcases}.

\section*{Acknowledgments}
The research of S. Du and S.N. Stechmann is partially supported
by grant NSF DMS 2324368
and
by the Office of the Vice Chancellor for Research and Graduate Education at the University of Wisconsin-Madison with funding from the Wisconsin Alumni Research Foundation.
{\color{black}The authors also thank the referees for their suggestions which have led to an improvement of the presentation of the paper.}

\appendix
\section{Implementation of HDG: further details}
\label{sec:app_hdg_mat}
Here we explain how the matrices in \eqref{eq:hdg_local_all_sub_mat} are assembled in more detail.
For the first term, we have
\begin{align*}
    (B^K[u_h])_{K,K^a,i} &=
    \sum_{F\in\mc E_K}\sum_{K_*^a\cdot\mb n_F\ge0}\int_{K_*^a}\int_F (\mb s\cdot\mb n_F)u_h v\\
    &=\sum_{F\in\mc E_K}\sum_{K_*^a\cdot\mb n_F\ge0}\int_{K_*^a}\int_F (\mb s\cdot\mb n_F)\sum_{j}u_{j}^{K\times K_*^a}\, \varphi_{j}^{K\times K_*^a} \varphi_{i}^{K\times K^a}\\
    &=\sum_{j}u_{j}^{K\times K^a}
    \left(
    \sum_{F\in\mc E_K}
    \mathbbm 1_{K^a\cdot\mb n_F\ge0}\int_{K^a}\int_F (\mb s\cdot\mb n_F)\, \varphi_{j}^{K\times K^a} \varphi_{i}^{K\times K^a}
    \right).
\end{align*}
For the second term, we have
\begin{align*}
    (\widehat B^K[\widehat u_h])_{K,K^a,i}&=\sum_{F\in\mc E_K}\sum_{K_*^a\cdot\mb n_F\le0}\int_{K_*^a}\int_F (\mb s\cdot\mb n_F)\widehat{u}_h v\\
    &=\sum_{F\in\mc E_K}\sum_{K_*^a\cdot\mb n_F\le0}\int_{K_*^a}\int_F (\mb s\cdot\mb n_F) \sum_{j}\widehat u_{j}^{F\times K_*^a}\, \psi_{j}^{F\times K_*^a}\varphi_i^{K\times K^a}\\
    &=\sum_{F\in\mc E_K}\sum_{j}\widehat u_{j}^{F\times K^a}\left(\mathbbm 1_{K^a\cdot\mb n_F\le0}\int_{K^a}\int_F (\mb s\cdot\mb n_F)  \psi_{j}^{F\times K^a}\varphi_i^{K\times K^a}\right).
\end{align*}
For the third term, we have
\begin{align*}
    (C^K[u_h])_{K,K^a,i}&
    =\sum_{K^a\in\mc T_h^a}\int_{K^a}\int_K u_h(\mb s\cdot\nabla v)\\
    &=\sum_{j}u_{j}^{K\times K^a}\left(\int_{K^a}   \int_K  \varphi_{j}^{K\times K^a}\mb s\cdot\nabla \varphi_i^{K\times K^a}\right).
\end{align*}
For the fourth term, we have
\begin{align*}
    (M^K[u_h])_{K,K^a,i}&=\sum_{K^a\in\mc T_h^a}\int_{K^a}\int_K\sigma_eu_hv\\
    &=\sum_ju_j^{K\times K^a}\left(\int_{K^a}\int_K\sigma_e\varphi_j^{K\times K^a}\varphi_i^{K\times K^a}\right).
\end{align*}
For the fifth term for scattering, we have
\begin{align*}
    (S^K[u_h])_{K,K^a,i}&=
    \sum_{K^a\in\mc T_h^a}\int_{K^a}\int_K\sigma_s(\mb x)\int_Sp(\mb s,\mb s')u_h(\mb x,\mb s')ds'v_h(\mb x,\mb s)dxds\\
    &=\int_{K^a}\int_K\sigma_s(\mb x)\int_Sp(\mb s,\mb s')
    \sum_{K_{*}^a}\sum_{j}u_{j}^{K\times K_{*}^a}\, \varphi_{j}^{K\times K_{*}^a}(\mb x,\mb s')ds'\varphi_i^{K\times K^a}(\mb x,\mb s)dxds\\
    &=\sum_{K_{*}^a\in\mc T_h^a}\sum_{j}u_{j}^{K\times K_{*}^a}
    \left(\int_K\sigma_s(\mb x)\int_{K^a}\int_{K_*^a}p(\mb s,\mb s')
    \varphi_{j}^{K\times K_{*}^a}(\mb x,\mb s')
    \varphi_i^{K\times K^a}(\mb x,\mb s)ds'dsdx\right).
\end{align*}

%\bibliographystyle{plain}
% \bibliographystyle{abbrv}
% \bibliography{ref}

\end{document}